\documentclass[reqno]{scrartcl}

\usepackage{egen,egenengelsk}
\usepackage[bookmarks]{hyperref}

\newcommand{\tsn}{T_{X^{[n]}}^{[\![n]\!]}}
\newcommand{\wh}{\widehat}

\newcommand{\CM}{CM}
\newcommand{\csm}{c_{\text{SM}}}
\newcommand{\relt}{T_{\ol{TX}^{[n]}}}
\newcommand{\dbl}{[\![}
\newcommand{\dbr}{]\!]}

\let\originalleft\left
\let\originalright\right
\renewcommand{\left}{\mathopen{}\mathclose\bgroup\originalleft}
\renewcommand{\right}{\aftergroup\egroup\originalright}

\title{Universal polynomials for tautological integrals on Hilbert schemes}
\author{Jørgen Vold Rennemo}
\date{}
\begin{document}
\maketitle
\begin{abstract}
We show that tautological integrals on Hilbert schemes of points can be written in terms of universal polynomials in Chern numbers. The results hold in all dimensions, though they strengthen known results even for surfaces by allowing integrals over arbitrary ``geometric'' subsets (and their Chern--Schwartz--MacPherson classes).

We apply this to enumerative questions, proving a generalised Göttsche Conjecture for all singularity types and in all dimensions. So if $L$ is a sufficiently ample line bundle on a smooth variety $X$, in a general subsystem $\PP^d \subset |L|$ of appropriate dimension the number of hypersurfaces with given singularity types is a polynomial in the Chern numbers of $(X,L)$. 

When $X$ is a surface, we get similar results for the locus of curves with fixed ``BPS spectrum'' in the sense of stable pairs theory.
\end{abstract}

\section{Results}
Let $X$ be a smooth, projective, connected variety of dimension $d$, and let $E$ be an algebraic vector bundle on $X$. Denote by $X^{[n]}$ the Hilbert scheme of length $n$ subschemes of $X$, and let $E^{[n]}$ be the tautological bundle on $X^{[n]}$ with fibre $H^0(Z, E|_Z)$ at $Z \in X^{[n]}$.

We study integrals of products of Chern classes of $E^{[n]}$ over what we call \emph{geometric subsets} of $X^{[n]}$. Geometric subsets form a natural class of subsets definable without reference to the global geometry of $X^{[n]}$. Specifically, the geometric subsets of $X^{[n]}$ are those generated by finite unions, intersections and complements from sets of the kind
\[
 \{Z \in X^{[n]} \mid Z = Z_1 \sqcup \ldots \sqcup Z_k,\ \ Z_i\text{ is of type }Q_i\}.
\]
Here each $Q_i$ is a constructible subset of the punctual Hilbert scheme $\Hilb^{n_i}_0(\CC^d)\subset \Hilb^{n_i}(\CC^d)$ of subschemes supported at $0\in \CC^d$, where we have $\sum n_i = n$. By ``$Z_i$ is of type $Q_i$'' we mean that $Z_i$ is isomorphic as an abstract scheme to some element of $Q_i$, so we require that the $i$-th connected component of $Z$ is of isomorphism type in a specified family $Q_i$.\footnote{Clearly, a geometric subset $P \subseteq X^{[n]}$ is constructible and has the property that if $Z \in P$, $Z^\pr \in X^{[n]}$ and $Z \cong Z^\pr$, then $Z^\pr \in P$. The converse is not quite true. For the basic subsets generating the algebra of geometric subsets we instead impose a similar condition on each connected component of $Z$, and this is a stronger requirement than the one above.

As an example of a subset with the above property which is not geometric by our definition, consider the following set. Let $X$ be a surface, and let $P \subset X^{[n]}$ be the set containing all $Z = Z_1 \sqcup Z_2$ such that (1) each $Z_i$ is defined by an ideal $(\fr{m}^5_i,f_i)$ where $f_i$ is a product of 4 distinct linear factors, and (2) the cross ratio of the factors of $f_1$ equals that of the factors of $f_2$.}

A $k$-variable Chern polynomial is a polynomial in the formal variables $c_i^{\left(j\right)}$, where $i \ge 1$ and $1 \le j \le k$. We treat such a Chern polynomial as a function from $k$-tuples of vector bundles to cohomology by the rule
\[
 c_i^{\left(j\right)}(E_1, \ldots, E_k) = c_i(E_j),
\]
extended linearly and multiplicatively to all Chern polynomials.

A Chern monomial is a monomial in the variables $c_i^{\left(j\right)}$. The \emph{weight} of a Chern monomial $c_{i_1}^{\left(j_1\right)}\cdots c_{i_k}^{\left(j_k\right)}$ is defined to be $\sum_{l=1}^k i_l$, so that treating a Chern monomial of weight $l$ as a function, its image will be in $H^{2l}(X)$. Denote by $\CM(k,l)$ the set of $k$-variable Chern monomials of weight $l$.

If $Y$ is a proper scheme and $Z \subset Y$ a constructible subset, we denote by $\csm(Z)\in H_*(Y)$ the Chern--Schwartz--MacPherson class of $Z$. The construction and basic properties of this class are reviewed in Section \ref{section:preliminaries}.

\begin{nthm}
\label{mainthm}
Let $X$ be a smooth, projective, connected variety of dimension $d$, $E$ an algebraic vector bundle on $X$, and $F$ a (1-variable) Chern polynomial. Let either
\begin{enumerate}
 \item[\textit{(i)}] $N = \deg (F(E) \cap [P])$, for $P\subseteq X^{[n]}$ closed and geometric, or
 \item[\textit{(ii)}] $N = \deg (F(E) \cap \csm(P))$, for $P\subseteq X^{[n]}$ geometric.
\end{enumerate}
Then, there exists a polynomial $G$ in the variables $\{x_M\}_{M\in \CM\left(2,d\right)}$, depending only on $F$, the rank of $E$ and the type of $P$, such that if we assign to $x_M$ the Chern number $\deg M(T_X,E) \cap [X]$, we have $N=G((x_M))$.

Moreover, if every point $Z \in P$ represents a subscheme with support in at most $m$ points, the degree of $G$ is at most $m$.
\end{nthm}

If $X$ is a surface or a curve, so that the Hilbert scheme is nonsingular, we can extend this as follows.
\begin{nthm}
\label{surfacethm}
Assume $\dim X = 1$ or $2$, and let $F$ be a 2-variable Chern polynomial. Theorem \ref{mainthm} then holds with $F(E)$ replaced by $F(T_{X^{[n]}}, E)$ everywhere.
\end{nthm}

An outline of the proof of the main theorem is given in Section \ref{outline}, and the formal proof occupies Sections \ref{section:proof} and \ref{section:lemmaproof}.

In Section \ref{section:generatingfun}, we show that a certain generating function for some Chern integrals of Theorems \ref{mainthm} and \ref{surfacethm} (i) can be given a particular product form.

The strategy of the proof of the main theorem is motivated by J. Li's paper \cite{li}, where he shows that the degree of the virtual fundamental class on the Hilbert scheme of points on a threefold $X$ is given by a universal polynomial in the Chern numbers of $X$. The outer structure of that proof, i.e.\ using the scheme $X^{[\![n]\!]}$ and approximating by classes defined via $X^{[\![\alpha]\!]}$, is well suited to our problem. Dealing with geometric subsets, the tautological bundles $E^{[n]}$ and the Chern--Schwartz--MacPherson class requires new ingredients.

A special case of Theorem \ref{surfacethm} has been proved by Ellingsrud, Göttsche and Lehn using a completely different method, see \cite{EGL}. In our terminology, they treat the case where $X$ is a surface and the geometric subset $P$ is the whole of $X^{[n]}$.

We note that the method of \cite{EGL} yields a recursion which computes the universal polynomial explicitly. In contrast, our method is nonconstructive and relies at a crucial point on the fact that an element in the cohomology ring of a Grassmannian is a polynomial in the Chern classes of the universal bundle. 
Lacking a method of obtaining information about this polynomial, there is no apparent way of turning our proof into an algorithm.

\subsection{Enumerative applications}
\subsubsection{Counting singular curves in surfaces}
The main motivation for our result is to generalise the result known as the Göttsche Conjecture, which by now has several proofs, see \cite{kazarian}, \cite{KST}, \cite{liu}, and \cite{tzeng}. We recall the statement of the conjecture. Fix a surface $S$ with a line bundle $L$ which is ``sufficiently ample'', e.g. $L = M^{\otimes N}$, where $M$ is a very ample line bundle and $N$ is a sufficiently large integer. The precise definition of sufficiently ample uses the concept of $N$-very ampleness, see Section \ref{applications}.

Let $\delta$ be a positive integer, and call a curve $\delta$-nodal if it has $\delta$ nodes and no other singularities. If $L$ is sufficiently ample, the locus of $\delta$-nodal curves in $|L|$ has the expected codimension $\delta$, so that in a general linear subsystem $\PP^\delta \subset |L|$ there is a finite number of $\delta$-nodal curves. The simplest form of the conjecture is then that there exists a  degree $\delta$ polynomial $N_\delta$ in 4 variables, independent of $S$ and $L$, such that the number of $\delta$-nodal curves equals
\[
 N_\delta\left(c_1\left(L\right)^2, c_1\left(L\right)c_1\left(S\right), c_1\left(S\right)^2,c_2\left(S\right)\right).
\]

Our main application is the generalisation of this result to the case of curves with more general specified singularity types. Our approach follows the idea of Göttsche used in \cite[Sec.\ 5]{gottsche} to reduce the problem of counting nodal curves to an integral on the Hilbert scheme. He defines a closed subset $W \subseteq S^{[3\delta]}$ and shows that the number of $\delta$-nodal curves in the linear system $\PP^\delta$ equals the degree of
\[
 c_{2\delta}\left(L^{[3\delta]}\right)\cap[W],
\]
assuming $L$ is $(5\delta - 1)$-very ample. This idea was used by Tzeng in her proof of the Göttsche Conjecture \cite{tzeng}, which uses degenerations of $S$ to show that the degree of the above class is a polynomial in the Chern numbers of $(S,L)$.

The set $W$ appearing above is geometric, hence our theorem yields a different proof of Tzeng's result. Since our main theorem deals with more general loci in the Hilbert scheme of points, we may generalise the statement of Tzeng's theorem, replacing $\delta$-nodal curves with curves having other specified singularity types.

\newcommand{\curveCountProp}{Let $S$ be a smooth, projective, connected surface, let $L$ be a line bundle on $S$, and let $T_1, \ldots, T_k$ be analytic singularity types. There are expected codimensions $d_i$ associated with each $T_i$, and we let $d = \sum d_i$.

There is an integer $N$ and a rational polynomial $G$ of degree $k$ in 4 variables, depending only on the $T_i$, such that if $L$ is $N$-very ample, then in a general $\PP^d \subseteq |L|$ the number of curves having precisely $k$ singularities of types $T_i$ is
\[
 G(c_1^2(L), c_1(L)c_1(S), c_2(S), c_1^2(S)).
\]
The same statement holds when the $T_i$ are topological rather than analytic singularity types.}
\begin{prop}[\ref{curveCountProp}]
\curveCountProp
\end{prop}

For the original problem of counting nodal curves, the numbers of curves having $k$ nodes form a generating function
\[
N(S,L) = \sum_{\delta \ge 0} N_\delta\left(c_1\left(L\right)^2, c_1\left(L\right)c_1\left(S\right), c_1\left(S\right)^2,c_2\left(S\right)\right) q^\delta,
\]
which was conjectured by Göttsche \cite[Prop.\ 2.3]{gottsche} and shown by Tzeng \cite[Thm.\ 1.3]{tzeng} to have a specific product form
\[
N(S,L) = B_1^{c_1^2(L)}B_2^{c_1(L)c_1(S)}B_3^{c_1^2(S)}B_4^{c_2(S)},\ \ \ B_i \in \QQ[[q]].
\]

We generalise this statement as Corollary \ref{cor:curveGen}: Fixing distinct types $T_i$, collect the universal polynomials for the number of curves having $m_i$ singularities of type $T_i$ in a generating function; this then admits a product expansion similar to the above.

Both Proposition \ref{curveCountProp} and Corollary \ref{cor:curveGen} have recently been obtained independently by Li and Tzeng via a generalisation of Tzeng's degeneration approach \cite{LT}.

\subsubsection{Counting singular hypersurfaces}
By the same method we are able to count singular hypersurfaces in arbitrary dimensions.
\newcommand{\divCountProp}{Let $X$ be a smooth, projective, connected variety, let $L$ be a line bundle on $X$, and let $T_1, \ldots, T_k$ be analytic singularity types. There are expected codimensions $d_i$ associated with each $T_i$, and we let $d = \sum d_i$.

There is an integer $N$ and a rational polynomial $G$ in the Chern numbers of $(X,L)$, depending only on the $T_i$, such that if $L$ is $N$-very ample, then in a general $\PP^d \subseteq |L|$ the number of divisors having precisely $k$ isolated singularities of types $T_i$ is given by $G$.}
\begin{prop}[\ref{divCountProp}]
\divCountProp
\end{prop}
As in the curve case, a generating function for these universal polynomials can be written in a product form similar to the one of Corollary \ref{cor:curveGen}.

\subsubsection{Counting curves with given BPS spectra}
A different application of the main result concerns the locus of curves in a $\PP^k \subset |L|$ having given ``BPS spectrum''. For a reduced, complete, locally planar curve $C$ with arithmetic genus $g(C)$ and geometric genus $\ol{g}(C)$, we consider the generating function
\[
H_C(q) := \sum_{k=0}^\infty \chi\left(C^{[k]}\right)q^k
\]
Pandharipande and Thomas \cite{PT} show that there are $n_{i,C} \in \ZZ$ for $i = \ol{g}(C), \ldots, g(C)$ such that
\[
H_C(q) = \sum_{i=\ol{g}(C)}^{g(C)} n_{i,C} q^{g(C)-i}(1-q)^{2i-2}.
\]
If $C$ is smooth, we have $H_C(q) = (1-q)^{2g-2}$, so this result can be interpreted as saying that in general $H_C(q)$ decomposes as a sum of $n_{i,C}$ copies of $q^{g-i}H_{C_i}(q)$ where $C_i$ is smooth of genus $i$. We define $m_{i,C} = n_{g(C)-i, C}$, and it is then easy to check that the sequence of integers $(m_{i,C})_{i=0}^\infty$ depends only on the analytic types of the singularities of $C$. We refer to the sequence $(m_{i,C})$ as the \emph{BPS spectrum} of $C$.

Recent work of Maulik, settling a conjecture of Oblomkov and Shende, shows that the BPS spectrum of $C$ is explicitly determined by the Milnor numbers and HOMFLY polynomials of the links of the singularities of $C$ \cite{maulik_stable_2012, OS}. 
As a consequence, the BPS spectrum depends only on the topological types of the singularities of $C$.

 We show the following proposition.
\newcommand{\BPSProp}{Let $S$ be a smooth, projective, connected surface, let $L$ be a line bundle on $S$, and let $k \in \ZZ_{\ge 0}$. Let $m = (m_i)_{i=0}^\infty$ be a BPS spectrum, and denote by $|L|_m\subseteq |L|$ the locus of curves with BPS spectrum $m$. 

There is an integer $N$ and a rational polynomial $G$ in 4 variables, depending only on $k$ and $m$, such that if $L$ is $N$-very ample, then for a general $\PP^k \subset |L|$ we have
\[
 \chi(\PP^k \cap |L|_m) = G(c_1(L)^2, c_1(L)c_1(S), c_1(S)^2, c_2(S)).
\]}
\begin{prop}[\ref{prop:BPS}]
 \BPSProp
\end{prop}

This generalises the approach of Kool, Shende and Thomas' proof of the Göttsche Conjecture \cite{KST}, which involves the special case of the above proposition where $m$ is the spectrum of a $\delta$-nodal curve; that is $m = (m_i)$ with $m_i = \binom{\delta}{i}$. 

We note that in the proof of Proposition \ref{prop:BPS} it is essential to be able to take integrals over general geometric subsets of $S^{[n]}$. This is in contrast to the argument of \cite{KST}, where the integrals needed were taken over the whole of $S^{[n]}$, and so were already computed in \cite{EGL}. 

\subsection{Conventions}
Homology and cohomology is singular with coefficients in $\QQ$. By the degree of a class in $H_*\left(X\right)$ we mean its pushforward to $H_*\left(\pt\right) \cong \QQ$. In dealing with algebraic subsets of Hilbert schemes we always give these the reduced scheme structure, and few functors are represented. 

If $m$ is some number defined in terms of the data $X,\left(E_i\right),P,F$ of the theorem, we will use the shorthand ``$m$ is universal'' to mean that there exists a polynomial in the variables $x_M$ computing $m$, depending only on $F$ and the type of $P$ as in the main theorem.
\subsection{Acknowledgements}
I thank Martijn Kool, Ragni Piene, my supervisor Richard Thomas, and Yu-jong Tzeng for valuable discussions and comments on this paper. 
In particular, Piene pointed out to me the results of \cite{KP} used in Section \ref{subsubsection:topSing}.

\section{Preliminaries}
\label{section:preliminaries}
Let $X$ be a smooth, projective, connected variety of dimension $d$, and let $E$ be an algebraic vector bundle on $X$. We give the definition of the tautological bundle $E^{[n]}$ and recall the construction of the Chern--Mather and Chern--Schwartz--MacPherson (CSM) classes. 
In Section \ref{sec:HilbOrdered} we introduce the scheme $X^{[\![n]\!]}$ and in Section \ref{section:geometric} we discuss the notion of geometric subsets of $X^{[n]}$ and $X^{[\![n]\!]}$.

\subsection{The tautological bundle}
Denote by $\cZ \subset X^{[n]} \times X$ the universal subscheme over $X^{[n]}$, and let $p: \cZ \to X$ and $q: \cZ \to X^{[n]}$ be the projections. The \emph{tautological bundle} $E^{[n]}$ on $X^{[n]}$ is defined as
\[
 E^{[n]} = q_*(p^*(E)).
\]
The flatness of $q$ implies that $E^{[n]}$ is locally free, and we see that the fibre of $E^{[n]}$ at a point $Z \in X^{[n]}$ is the vector space $H^0\left(Z, E|_Z\right)$.

\subsection{Chern classes}
We next review the Chern--Mather and Chern--Schwartz--MacPherson classes. These classes are generalisations to singular varieties of the Poincaré dual of $c_\bullet(T_Y)$ for a smooth, proper $Y$, so for such $Y$ we have
\[
 \csm(Y) = c_M(Y) = c_\bullet(T_Y) \cap [Y].
\]

\subsubsection{Chern--Mather class}
While the Chern--Mather class does not appear in the main theorem, it is used in the definition of the CSM class and throughout the proof of the main theorem.

Let $Y$ be a be a reduced and irreducible projective scheme.
The first step is to construct the Nash blow up $\wt{Y} \to Y$. Suppose for a moment that $Y$ is \emph{affine}, reduced and irreducible of dimension $d$. Fix an embedding $f: Y \to \AA^N$, and let $Y_{\text{ns}}$ be the nonsingular part of $Y$. The tangent map $T_{Y_{ns}} \to f^*(T_{\AA^N})$ induces a morphism $g : Y_{ns} \to \Gr(d,N)$, and we take $\wt{Y}$ to be the closure of the graph $\Gamma_g \subset Y \times \Gr\left(d,N\right)$. The morphism $\wt{Y}\to Y$ is defined by the projection $Y\times \Gr\left(d,N\right) \to Y$, and we define the rank $d$ vector bundle $T_{\wt{Y}}$ on $\wt{Y}$ by restricting the universal bundle on $\Gr\left(d,N\right)$.

It can be shown that this construction is independent of the choice of affine embedding and globalises so that for any reduced, equidimensional scheme $Y$ we get a well defined $Y$-scheme $\wt{Y}$ with a bundle $T_{\wt{Y}}$. The morphism $\wt{Y} \to Y$ is the \emph{Nash blow up} of $Y$ and the bundle $T_{\wt{Y}}$ is the \emph{Nash bundle}.

\begin{ndefn}
The Chern--Mather class $c_M\left(Y\right) \in H_*\left(Y\right)$ is the pushforward of $c_\bullet\left(T_{\wt{Y}}\right) \cap \left[\wt{Y}\right]$ along $\wt{Y} \to Y$.
\end{ndefn}

\subsubsection{Chern--Schwartz--MacPherson class}
We recall the definition and basic properties of the Chern--Schwartz--MacPherson class. For details see \cite[Ex. 19.1.7]{fulton} and \cite{mac}.

Let $Y$ be a projective scheme, let $Z_*(Y)$ be the group of all cycles on $Y$, and let $F_*\left(Y\right)$ denote the group of constructible functions, where a function $f: Y \to \ZZ$ is called constructible if there exists a finite partition of $Y$ into constructible sets such that $f$ is constant on each set of the partition. Given any reduced scheme $V$, the local Euler obstruction $\nu_V :V \to \ZZ$ is a canonical constructible function determined at a point $x \in V$ by the local analytic structure of $V$ at $x$.\footnote{If $p: \wt{V} \to V$ is the Nash blow-up, we have $\nu_V(x) = \deg\left(c_\bullet\left(T_{\wt{V}}|_{p^{-1}(x)}\right) \cap s(p^{-1}(x),\wt{V})\right)$, where $p^{-1}(x)$ is the scheme-theoretic inverse image and $s(-,-)$ denotes the Segre class.}

There is a group homomorphism $\Omega : Z_*\left(Y\right) \to F_*\left(Y\right)$ defined on a primitive cycle $V$ by
\[
\Omega(V)(x) = 
\begin{cases}
\nu_V(x)\text{ if }x \in V \\
0\text{ if }x\not\in V
\end{cases}
\]
The map $\Omega$ is an isomorphism, as can be shown using the fact that $\nu_V(x) = 1$ if $x$ is a nonsingular point of $V$. The Chern--Mather class defines a homomorphism $c : Z_*\left(Y\right) \to H_*\left(Y\right)$ by letting
\[
c\left(V\right) = i_*\left(c_M\left(V\right)\right),
\]
where $i: V \to Y$ is the inclusion of a primitive cycle. 

The Chern--Schwartz--MacPherson class $\csm(f)$ is now defined for any constructible function $f$ by
\[
\csm(f) = c(\Omega^{-1}(f)).
\]
It is clear that $\csm$ is a homomorphism. If $Z \subset Y$ is a constructible subset, we write $\csm(Z) = \csm(1_Z)$, where $1_Z$ is the characteristic function of $Z$. Additivity of $\csm$ translates to
\[
 \csm(Z_1 \cup Z_2) = \csm(Z_1) + \csm(Z_2) - \csm(Z_1 \cap Z_2).
\]

Given a morphism of proper schemes $g: Y_1 \to Y_2$, one can define a homomorphism $g_*:F(Y_1) \to F(Y_2)$ by letting
\[
 g_*(1_V)(x) = \chi(g^{-1}(x) \cap [V])\ \ \ \ x \in Y_2,
\]
where $V\subset Y_1$ is a primitive cycle, and $\chi$ is the topological Euler characteristic.\footnote{The fact that this is well defined is shown in \cite{mac}.} The main property of CSM classes, shown in \cite{mac}, is that if $g$ is proper, we have $g_*(\csm(f)) = \csm(g_*(f))$. Appying this to the case where $Y$ is proper and $g$ is the map $Y \to \pt$, we find that $\deg \csm(Y) = \chi(Y)$.

\subsection{The Hilbert scheme of ordered points}
\label{sec:HilbOrdered}
Following \cite{li}, we introduce the scheme $X^{[\![n]\!]}$, which will play an essential role in the proof of Theorems \ref{mainthm}/\ref{surfacethm}.
\begin{ndefn}
\label{def:orderedpoints}
The \emph{Hilbert scheme of ordered points}, denoted $X^{[\![n]\!]}$, is the scheme defined by the Cartesian diagram
\[
 \begin{CD}
  X^{[\![n]\!]} @>>> X^{[n]} \\
  @VVV           @VVV    \\
  X^n       @>>> \Sym^n(X),
 \end{CD}
\]
where the right hand arrow is the Hilbert-Chow morphism taking a subscheme $Z$ to its support cycle. We denote by $(Z,(x_i))$ the point in $X^{[\![n]\!]}$ mapping to $Z \in X^{[n]}$ and $(x_i) \in X^n$.
\end{ndefn}
Since $X^{[\![n]\!]} \to X^{[n]}$ is finite, we can reduce questions about the degree of homology classes on $X^{[n]}$ to questions about similar classes on $X^{[\![n]\!]}$.
Roughly speaking, the advantage of introducing $X^{[\![n]\!]}$ is that it naturally maps to $X^n$.
This makes it easier to handle than $X^{[n]}$, which maps to the more complicated scheme $\Sym^n(X)$. 

\subsection{Geometric subsets}
\label{section:geometric}
We now give the definition of geometric subsets of $X^{[n]}$ and of $X^{[\![n]\!]}$, along with some results on these which will be needed later.

Let $\Hilb_0^n(\CC^d)$ be the punctual Hilbert scheme, defined as the closed subset of $\Hilb^n(\CC^d)$ parametrising subschemes supported at the origin. We define punctual geometric subsets to be the constructible subsets of the punctual Hilbert scheme containing all 0-dimensional schemes of given isomorphism types.
\begin{ndefn}
A punctual geometric set is a subset $Q \subseteq \Hilb_0^n\left(\CC^d\right)$ which is constructible and satisfies: If $Z \in Q$ and $Z^\pr \in \Hilb^n_0(\CC^d)$ are such that $Z \cong Z^\pr$ (as abstract $\CC$-schemes), then $Z^\pr \in Q$.
\end{ndefn}

A collection of punctual geometric subsets will naturally defines a subset of $X^{[n]}$:
\begin{ndefn}
\label{def:standardgeometric} Let $Q_1, \ldots, Q_k$ be punctual geometric sets such that $Q_i \subseteq \Hilb_0^{n_i}(\CC^d)$, and let $n = \sum n_i$.

We define $P(Q_1, \ldots, Q_k) \subseteq X^{[n]}$ to be the set of all $Z = Z_1 \sqcup \ldots \sqcup Z_k$, where every $Z_i$ is isomorphic to a $Z_i^\pr \in Q_i$.
\end{ndefn}

If we additionally specify how to label these $Z$, we obtain a subset of $X^{[\![n]\!]}$:
\begin{ndefn}
\label{def:standardLabelledGeometric}
Let $Q_i$ and $n_i$ be as above. Let $A = (A_1, \ldots ,A_k)$ be a $k$-tuple of subsets of $\{1, \ldots, n\}$, such that $|A_i| = n_i$, and such that the $A_i$ define a partition of $\{1, \ldots, n\}$.

We define $T(Q_1, \ldots, Q_k; A) \subseteq X^{[\![n]\!]}$ to be the set of all $(Z, (x_i))$ such that $Z = Z_1 \sqcup \ldots \sqcup Z_k$, where every $Z_i$ is isomorphic to a $Z_i^\pr \in Q_i$, and such that $x_i = \Supp Z_j$ if $i \in A_j$.
\end{ndefn}

We can now give the definition of geometric subsets of $X^{[n]}$ and $X^{[\![n]\!]}$:
\begin{ndefn}\mbox{}
\begin{itemize}
\item A subset $P \subseteq X^{[n]}$ is geometric if it can be expressed using sets of the form $P(Q_1, \ldots, Q_k)$ and a finite composition of the operations union, intersection and complement.
\item A subset $T \subseteq X^{[\![n]\!]}$ is geometric if it can be expressed using sets of the form $T(Q_1, \ldots, Q_k; A)$ and a finite composition of the operations union, intersection and complement.
\end{itemize}
\label{def:geometric}
\end{ndefn}
An equivalent definition which will be convenient is the following.
\begin{ndefn}\mbox{}
\begin{itemize}
\item A subset $P \subseteq X^{[n]}$ is geometric if it can be expressed using sets of the form $\ol{P(Q_1, \ldots, Q_k)}$, where the $Q_i$ are closed and irreducible, together with a finite composition of the operations union, intersection and complement.
\item A subset $T \subseteq X^{[\![n]\!]}$ is geometric if it can be expressed using sets of the form $\ol{T(Q_1, \ldots, Q_k; A)}$, where the $Q_i$ are closed and irreducible, together with a finite composition of the operations union, intersection and complement.
\end{itemize}
\label{def:geometricalt}
\end{ndefn}
The equivalence of Definitions \ref{def:geometric} and \ref{def:geometricalt} is shown in Lemma \ref{lemma:geometric} (vii).

\begin{ex}
The only geometric subsets of $X^{[1]} = X$ are $\emptyset$ and $X$. In $X^{[2]}$ there are 4 geometric subsets: The sets $\emptyset$, $X^{[2]}$, the set parametrising pairs of disjoint points, and the set parametrising length 2 subschemes with support in one point. 

When $X$ is a surface, a naturally occurring example of a geometric subset is the set $W \subset X^{[3\delta]}$, defined as the closure of
\[
 \{Z \in X^{[3\delta]} \mid Z = Z_1 \sqcup \cdots \sqcup Z_\delta,\ \ Z_i = \Spec \cO_{S,x_i}/\mathfrak{m}_{x_i}^2\}.
\]
This set appears in Tzeng's proof of the Göttsche Conjecture.
\end{ex}

Clearly, if $P \subseteq X^{[n]}$ is geometric and $Z \in P$, then for any $Z^\pr \in X^{[n]}$ such that $Z \cong Z^\pr$ we have $Z^\pr \in P$. In other words, a geometric subset is a union of isomorphism classes of subschemes $Z \in X^{[n]}$.

We have a notion of isomorphism between points of $X^{[\![n]\!]}$, defined by saying that $(Z,(x_i)) \cong (Z^{\pr},(x_i^\pr))$ if there exists an isomorphism $Z \cong Z^\pr$ which takes $Z|_{x_i}$ to $Z^\pr|_{x_i^\pr}$ for every $i$. Then similarly a geometric subset of $X^{[\![n]\!]}$ is a union of isomorphism classes of pairs.

Let $X$ and $Y$ be smooth varieties of equal dimension, and let $P$ and $Q$ be geometric subsets of $X^{[n]}$ and $Y^{[n]}$ respectively. We say that $P$ and $Q$ are of the same \emph{type} if the isomorphism classes of the points in $P$ are the same as the isomorphism classes of points in $Q$. Clearly, for any geometric subset of $X^{[n]}$ there is a unique geometric subset of $Y^{[n]}$ of the same type. The type of a geometric subset $T \subseteq X^{[\![n]\!]}$ is defined in the same way.

The following lemma contains some elementary facts about geometric subsets.
\begin{nlemma} 
\label{lemma:geometric}
Let $P \subseteq X^{[n]}$ and $T \subseteq X^{[\![n]\!]}$ as sets, and let $p : X^{[\![n]\!]} \to X^{[n]}$ be the natural forgetful morphism. 
 \begin{enumerate}
  \item[(i)] $P$ is geometric $\Leftrightarrow$ $p^{-1}(P)$ is geometric.
  \item[(ii)] $T$ is geometric $\Rightarrow$ $p(T)$ is geometric.
  \item[(iii)] $P$ is geometric $\Leftrightarrow$ $P$ is a finite union of sets of the form $P((Q_i))$.
  \item[(iv)] $T$ is geometric $\Leftrightarrow$ $T$ is a finite union of sets of the form $T((Q_i);A)$.
  \item[(v)] $P$ is geometric, closed and irreducible $\Leftrightarrow$ $P$ is of the form $\ol{P((Q_i))}$ for closed, irreducible $Q_i$.
  \item[(vi)] $T$ is geometric, closed and irreducible $\Leftrightarrow$ $T$ is of the form $\ol{T((Q_i);A)}$ for closed, irreducible $Q_i$.
  \item[(vii)] Definitions \ref{def:geometric} and \ref{def:geometricalt} are equivalent.
 \end{enumerate}
\end{nlemma}
\begin{proof}
In this proof, ``geometric subsets'' means a set satisfying Definition \ref{def:geometric}.

(iv) It suffices to show that intersections and complements of sets of the form $T((Q_i);A)$ are expressible as unions of such sets. Let $T((Q_i);A)$ and $T((Q_i^\pr);A^\pr))$ be geometric sets with $A = (A_1, \ldots, A_k)$ and $A^\pr = (A_1^\pr, \ldots A_l^\pr)$. Then, if $T((Q_i);A) \cap T((Q_i^\pr);A^\pr) \not= \emptyset$, we have $k=l$ and the $k$-tuple $A$ is a permutation of the $k$-tuple $A^\pr$. In this case, we may relabel the indices of the $A^\pr_i$ to get $A = A^\pr$, and then $T((Q_i);A) \cap T((Q_i^\pr);A)) = T((Q_i\cap Q_i^\pr);A)$. 

Next we see that for any $T((Q_i);A)$ the set $X^{[\![n]\!]} \sm T((Q_i);A)$ is the union of all sets $T((\Hilb_0^{n_i}(\CC^d));A^\pr)$ where $A^\pr$ is not a permutation of $A$, and the sets 
\[
T(\Hilb^{n_1}_0(\CC^d), \ldots, \Hilb^{n_i}_0(\CC^d)\sm Q_i, \ldots ,\Hilb^{n_k}_0(\CC^d); A) 
\]
for $i = 1, \ldots, k$.

(ii) Follows from (iv) and $p(T((Q_i);A)) = P((Q_i))$.

(i) The $\Rightarrow$ follows from the fact that $p^{-1}(P(Q_i))$ is the union of $T((Q_i);A)$ for all admissible $A$. The $\Leftarrow$ implication then follows from (ii) and surjectivity of $p$.

(iii) Follows from (i), (ii), (iv) and the surjectivity of $p$.

(v) By (iii), we may write $P = \cup_j P((Q_{i,j})_i)$, and as $P$ is closed, we have $P = \cup_j \ol{P((Q_{i,j})_i)} = \cup_j \ol{P((\ol{Q_{i,j}})_i)}$. Irreducibility of $P$ implies $P = \ol{P((\ol{Q_{i,j}})_i)}$ for some $j$, so we may take $Q_i = \ol{Q_{i,j}}$. It remains to show that the $Q_i$ can be chosen to be irreducible. Suppose not, then we have for instance $Q_1$ reducible. Let $Q_1 = \cup_j Q_{1,j}$ be the decomposition of $Q_1$ into closed, irreducible subsets. Each $Q_{1,j}$ must be equal to the closure of its orbit under the natural action of $\Aut(\cO_{\AA^d,0}/\fr{m}_0^{n_1})$ on $\Hilb^{n_1}_0(\CC^d)$, hence we see that the $Q_{1,j}$ are geometric. 

We then have $P = \cup_j \ol{P(Q_{1,j}, Q_2, \ldots, Q_k)}$, and as $P$ is irreducible we may replace $Q_1$ with some $Q_{1,j}$. Repeat to get all $Q_i$ irreducible, proving the $\Rightarrow$ implication. The $\Leftarrow$ implication is easy and omitted, but note that it depends on the hypothesis that $X$ is connected.

(vi) Similar to (v).

(vii) It is obvious that a $P$ satisfying Def.\ \ref{def:geometricalt} satisfies Def.\ \ref{def:geometric}. For the converse, note that the closed, geometric $P$ generate all geometric subsets by unions, intersections and complements. The proof of (v) shows that a closed, geometric $P$ is the union of sets of the form $\ol{P((Q_i))}$ with $Q_i$ closed and irreducible. Hence closed, geometric $P$ satisfy Definition \ref{def:geometricalt}, and the claim follows. The case of $T$ is similar.
\end{proof}

\section{Outline of proof}
\label{outline}
We give an outline of the proof of the main theorem. 
We restrict our attention in the outline to Theorem \ref{mainthm} (i), ignoring the extra complications of (ii) and Theorem \ref{surfacethm}. Arguing as in the proof of Lemma \ref{lemma:geometric} (v), we see that the irreducible components of $P$ are geometric of type depending only on the type of $P$. Hence we may assume that $P$ is irreducible.

The set-up is then that we are given a closed, irreducible, geometric subset $P$ of $X^{[n]}$, a Chern polynomial $F$ and a vector bundle $E$, and we want to show that
\[
 \deg F\left(E^{[n]}\right) \cap [P]
\]
is given by a universal polynomial.

\subsection{Reduction to $X^{[\![n]\!]}$}
The first step is to replace $X^{[n]}$ with the Hilbert scheme of ordered points $X^{[\![n]\!]}$. 

Define the bundle $E^{[\![n]\!]}$ on $X^{[\![n]\!]}$ as the pullback of $E^{[n]}$ along $X^{[\![n]\!]} \to X^{[n]}$. 
In Lemma \ref{covering} we construct a closed, irreducible $T \subseteq X^{[\![n]\!]}$ which is geometric, maps properly and finitely onto $P$, and which is such that $\deg\left(T/P\right)$ and the type of $T$ is determined by the type of $P$. The projection formula then gives
\[
\deg\left(T/P\right)\left(\deg F\left(E^{[n]}\right) \cap [P]\right) = \deg F\left(E^{[\![n]\!]}\right) \cap [T].
\]
Thus, it suffices to show that $\deg F\left(E^{[\![n]\!]}\right) \cap [T]$ is given by a universal polynomial.

\subsection{Approximating spaces}
\label{sec:ApproximatingSpaces}
Following \cite{li}, we let $\alpha$ be a partition of $\{1, \ldots, n\}$, and define the scheme $X^{[\![\alpha]\!]}$. Considering $\alpha$ as a set of subsets of $\{1, \ldots, n\}$, we let
\[
 X^{[\![\alpha]\!]} = \prod_{A \in \alpha} X^{[\![|A|]\!]}.
\]

So, for example, if $\alpha$ is the partition of $\{1, \ldots, n\}$ into $n$ one-element sets, we have $X^{[\![\alpha]\!]} = X^n$. 
At the other extreme, for the trivial partition $\Lambda = \{\{1, \ldots, n\}\}$, we have $X^{[\![\Lambda]\!]} = X^{[\![n]\!]}$. 
In general, the scheme $X^{[\![\alpha]\!]}$ parametrises ordered collections of $n$ points in $X$, with the additional data that when $k$ points with labels in the same set in the partition $\alpha$ come together at $x$, one must specify a length $k$ subscheme supported at $x$.

Consider the subset of $X^{[\![n]\!]}$ where no two points with labels in different sets of $\alpha$ come together. 
This subset is naturally isomorphic to a subset of $X^{[\![\alpha]\!]}$, so we get a rational map $f:X^{[\![n]\!]} \dashrightarrow X^{[\![\alpha]\!]}$.

We define a bundle $E^{[\![\alpha]\!]}$ on $X^{[\![\alpha]\!]}$ such that $f^*(E^{[\![\alpha]\!]}) = E^{[\![n]\!]}$ on the locus where $f$ is defined. If $\alpha$ is such that $T$ intersects the locus where $f$ is defined, we get a closed, irreducible subset $T_\alpha := \ol{f(T)} \subseteq X^{[\![\alpha]\!]}$. We have $E^{[\![\Lambda]\!]} = E^{[\![n]\!]}$ and $T_\Lambda = T$. 

All the $T_\alpha$ are birational to $T$.
We let $Q$ be the closure of the image of $T \dashrightarrow \prod_{\alpha} T_\alpha$, where the product is over all $\alpha$ such that $T_\alpha$ is defined. 
The projections from $\prod T_\alpha$ induce proper, birational morphisms from $Q$ to every $T_\alpha$.

\subsection{Approximating cohomology classes}
In what follows we restrict attention to partitions $\alpha$ such that $T_\alpha$ is defined. We define the class $C_\alpha \in H^*(Q)$ by
\[
 C_\alpha = F\left(E^{[\![\alpha]\!]}\right),
\]
suppressing the pullback of $E^{[\![\alpha]\!]}$ along $Q \to T_\alpha$. Let $C = C_\Lambda$. By the projection formula, $\deg \left(C \cap [Q]\right) = \deg \left(F\left(E^{[\![n]\!]}\right) \cap [T]\right)$, so the proof of the main theorem is reduced to showing that $\deg C \cap [Q]$ is universal.

We next define a class 
\begin{equation*}
\label{eq:DefOfD}
D = C + \sum_{\alpha \not= \Lambda} k_\alpha C_\alpha,
\end{equation*}
where the $k_\alpha$ are certain combinatorially defined integers (cf.\ \cite[Sec.\ 5.4]{li}). There is a natural morphism $Q \to X^n$, and one should think of the class $D$ as being supported on (a neighbourhood of) the set $Q|_\Delta$, where $\Delta \subset X^n$ is the small diagonal. The choice of the integers $k_\alpha$ is motivated by the fact that they make $D$ vanish on the complement of this locus.

For any $\alpha \not= \Lambda$, the scheme $X^{[\![\alpha]\!]}$ is by definition a product of schemes $X^{[\![m]\!]}$ with $m < n$. This induces product decompositions of $E^{[\![\alpha]\!]}$ and $T_\alpha$, which allow us to express $\deg \left(C_\alpha \cap [Q]\right)$ in terms of integrals of Chern classes of $E^{[\![m]\!]}$ over geometric subsets of $X^{[\![m]\!]}$ with $m<n$. By induction on $n$ we can thus show that $\deg \left(C_\alpha \cap [T_\alpha]\right)$ is universal for $\alpha \not =\Lambda$. This argument gives Lemma \ref{product}, by which it suffices to show that 
\[
\deg \left(D \cap [Q]\right)
\]
is universal.

\subsection{Relative constructions}
Consider the tangent bundle $TX \to X$, and let $\ol{TX} := \PP(\cO_X \oplus TX)$ be the natural compactification. Let $\Hilb^n(\ol{TX}/X)$ be the relative Hilbert scheme, which parametrises length $n$ subschemes of fibres of $\ol{TX} \to X$. Emulating the definition of $Q$ and $E^{[\![\alpha]\!]}$ with $\Hilb^n(\ol{TX}/X)$ replacing $X^{[n]}$, we define the scheme $\cQ$ and the bundles $\cE^{[\![\alpha]\!]}$ on $\cQ$. The classes $\cC_\alpha, \cD \in H^*(\cQ)$ are defined similarly to $C_\alpha$ and $D$.

Denote by $\ol{TX}^n$ the $n$-fold fibre product of $\ol{TX}$ over $X$. There are natural morphisms $g : Q \to X^n$ and $h:\cQ \to \ol{TX}^n$, where $g$ is given by composing $Q \to X^{[\![n]\!]}$ and $X^{[\![n]\!]} \to X^n$, and $h$ is defined similarly. Let $\Delta \subset X^n$ be the small diagonal, and consider $X$ as a subset of $TX^n \subset \ol{TX}^n$ using $n$ copies of the 0-section. 

Let $U\subset Q$ and $\cU \subset \cQ$ be Euclidean open neighbourhoods of $g^{-1}\left(\Delta\right)$ and $h^{-1}\left(X\right)$, respectively. Choosing $U$ and $\cU$ small enough, we can define a topological isomorphism $f : U \to \cU$, inducing topological isomorphisms of the bundles $f^*\left(\cE^{[\![\alpha]\!]}\right) \cong E^{[\![\alpha]\!]}$. 

To define the map $f$, we follow \cite{li}, and begin with an exponential map $\exp: TX \to X^2$. This map is defined in a neighbourhood of the 0-section $X \subset TX$, and maps a neighbourhood of the 0-section homeomorphically onto a neighbourhood of the diagonal $\Delta \subset X^2$. 
The map $\exp$ is analytic when restricted to a fibre of $TX$.

The inverse of this map induces a local homeomorphism $X^{[n]} \to \Hilb^n(\ol{TX}/X)$, defined in a neighbourhood of the locus of subschemes supported at a single point. 
Tracing through the parallel steps in the definitions of $Q$ and $\cQ$ gives a local homeomorphism $f: Q \to \cQ$. 
The details of this construction are contained in Lemmas \ref{f1}-\ref{f3}.

Recall the class $D = C + \sum_{\alpha \not= \Lambda} k_\alpha C_\alpha$. The bundles $E^{[\![\alpha]\!]}$ are canonically isomorphic over various open subsets of $Q$, and we use these isomorphisms to show that the terms of this sum cancel out when restricted to $Q \sm U$. In other words, $D$ vanishes when restricted to $Q \sm U$. Similarly $\cD$ vanishes upon restriction to $\cQ \sm \cU$.

Hence the classes $D$ and $\cD$ are concentrated over $U$ and $\cU$, respectively. In particular, there are relative cohomology classes $\ol{D} \in H^*(Q, Q\sm U)$ and $\ol{\cD} \in H^*(\cQ, \cQ \sm \cU)$ lifting $D$ and $\cD$.

There is a map $f^*: H^*(\cQ, \cQ\sm \cU) \to H^*(Q, Q\sm U)$, defined by excision, possibly after shrinking $U$ and $\cU$. Lemma \ref{rel} shows that we can choose $\ol{D}$ and $\ol{\cD}$ in such a way that $f^*(\ol{\cD}) = \ol{D}$. Hence we get
\[
 \deg D \cap [Q] = \deg \ol{D} \cap [Q, Q \sm U] = \deg \ol{\cD} \cap [\cQ, \cQ \sm \cU] = \deg \cD \cap [\cQ].
\]
The proof of Lemma \ref{rel} is quite technical and occupies Section \ref{section:lemmaproof}.

\subsection{Pullback from the Grassmannian}
The final step in the argument is to show that $\deg \cD \cap [\cQ]$ is universal. 

Let $U \to \Gr(d,N)$ be the universal rank $d$ subbundle over a Grassmannian. Here $N$ is any integer large enough that $TX$ embeds as a topological subbundle of $\cO_X^N$, so that there is a continuous classifying map $f: X \to \Gr(d,N)$ with $TX \cong f^*(U)$ as topological bundles. We define the scheme $\cQ_{\Gr}$ by the same construction as $\cQ$, replacing $\Hilb^n(\ol{TX}/X)$ with $\Hilb^n(\ol{U}/\Gr\left(d,N\right))$ throughout. There is a natural morphism $\cQ_{\Gr} \to \Gr(d,N)$ and the Cartesian diagram of topological spaces
\[
\begin{CD}
\cQ @>>>  \cQ_{\Gr} \\
@VpVV      @VVV \\
X   @>f>> \Gr(d,N).
\end{CD}
\]
Let $h: X \to \Gr(r,M)$ be a continuous classifying map for $E$, and consider the Cartesian diagram
\[
 \begin{CD}
  \cQ @>g>> \cQ_{\Gr} \times \Gr(r,M) \\
  @VpVV     @VVV    \\
  X   @>(f,h)>> \Gr\left(d,N\right) \times \Gr\left(r, M\right).
 \end{CD}
\]
We emphasise that the horizontal arrows in these diagrams are not required to be analytic. We define bundles $\cE_{\Gr}^{[\![\alpha]\!]}$ on $\cQ_{\Gr} \times \Gr(r,M)$ such that $g^*\left(\cE^{[\![\alpha]\!]}_{\Gr}\right) = \cE^{[\![\alpha]\!]}$. 

From the above it follows that there is a class $\cG \in H^*(\Gr\left(d,N\right) \times \Gr\left(r,M\right))$, depending only on $F$ and the type of $T$, such that
\[
 p_*(\cD \cap [\cQ]) = (f,h)^*(\cG) \cap [X].
\]
Now, as the rational cohomology ring of a Grassmannian is generated by the Chern classes of its universal bundle, $(f,h)^*(\cG)$ is a polynomial in the Chern classes of $T_X$ and $E$. We show that this polynomial is independent of the choice of $N$ and $M$. Hence $\deg \cD \cap [\cQ]$ is a universal linear combination of the Chern numbers of $(X,E)$, which concludes the proof of the main theorem.

\section{Proof of main theorem}
\label{section:proof}
We now begin the formal proof of the main theorem. 
To avoid dealing with Theorem \ref{mainthm} and \ref{surfacethm} separately, we adopt the following convention: When a formula includes $T_{X^{[n]}}$, terms involving $T_{X^{[n]}}$ should be ignored unless $\dim X = 1,2$, and so should all other statements involving $T_{X^{[n]}}$. 

For part (i), it suffices to treat the case where $P$ is irreducible. If $P$ is closed and irreducible, we have $\csm(P) = [P] + $ terms of lower homological degree, hence part (ii) implies part (i). By the inclusion-exclusion property of CSM classes, it finally suffices to prove (ii) for closed and irreducible $P$.

In this section and the next we have $X, P, E, d$ and $F$ as in the main theorem, and we assume that $P$ is closed and irreducible.

\subsection{Reduction to $X^{[\![n]\!]}$}
Recall the Hilbert scheme of ordered points $X^{[\![n]\!]}$, defined by the Cartesian diagram
\[
\begin{CD}
 X^{[\![n]\!]} @>>> X^{[n]} \\
 @VVV           @VVV \\
 X^n       @>>> \Sym^n(X),
\end{CD}
\]
where the right hand arrow is the Hilbert-Chow morphism.

\begin{ndefn}
Denote the pullbacks of $E^{[n]}$ and $T_{X^{[n]}}$ along $X^{[\![n]\!]}\to X^{[n]}$ by $E^{[\![n]\!]}$ and $T_{X^{[n]}}^{[\![n]\!]}$, respectively.
\end{ndefn}

We will use the projection formula to relate the degree of $F(E^{[n]}, T_{X^{[n]}}) \cap [P]$ to a similar class involving $E^{[\![n]\!]}$ and $\tsn$ on $X^{[\![n]\!]}$. The first step is to produce a closed, irreducible geometric subset $T \subseteq X^{[\![n]\!]}$ mapping finitely onto $P$ with universal degree.

\begin{nlemma}
\label{covering}
There exists a closed, irreducible geometric $T \subset X^{[\![n]\!]}$, such that $X^{[\![n]\!]} \to X^{[n]}$ maps $T$ finitely onto $P$. Up to a permutation of $\{1, \ldots, n\}$, the type of this $T$ is uniquely determined by the type of $P$.
\end{nlemma}
\begin{proof}
 As $P$ is closed and irreducible, by Lemma \ref{lemma:geometric} (v) we have $P = \ol{P((Q_i))}$ for closed and irreducible punctual geometric subsets $Q_i\subset\Hilb^{n_i}_0(\CC^d)$. We take $A = (A_1, \ldots, A_k)$ to be a $k$-element partition of $\{1, \ldots n\}$ such that $|A_k| = n_k$, and let $T = \ol{T((Q_i); A)}$. By Lemma \ref{lemma:geometric} (vi), we see that $T$ must have this form, hence the second claim follows.
\end{proof}

Recall that $c_M(T)$ denotes the Chern--Mather class of $T$.
\begin{nlemma}
In order to prove the main theorem, it suffices to prove the following statement: If $T$ is a closed, irreducible, geometric subset of $X^{[\![n]\!]}$, then
\[
 \deg F\left(E^{[\![n]\!]}, \tsn\right) \cap c_M(T)
\]
is given by a universal polynomial depending only on $F$ and the type of $T$, such that the degree of the polynomial is at most $l$, where $l$ is the maximum number of components of $Z$ for $(Z,(x_i)) \in T$.
\end{nlemma}
\begin{proof}
Let $p : X^{[\![n]\!]} \to X^{[n]}$ be the natural morphism. As explained at the beginning of this section, it is enough to show the main theorem in the case when $P \subseteq X^{[n]}$ is closed and irreducible. Let $T \subseteq X^{[\![n]\!]}$ be a closed, irreducible geometric set mapping finitely onto $P$, as provided by Lemma \ref{covering}.

We claim that the level sets of the local Euler obstruction $\nu_T$ are geometric. To see this, first write $T$ as a union of sets of the form $T((Q_i);A)$, as can be done by Lemma \ref{lemma:geometric} (iv). Consider a point $(Z, (x_j)) \in T$, and suppose that $(Z,(x_j)) \in T((Q_i); A) \subseteq T$ with $A = (A_1, \ldots, A_k)$. We have $Z = \sqcup_i Z_i$ with each $Z_i$ isomorphic to an element of $Q_i$ and $x_j = \Supp Z_i$ for $j \in A_i$. The local analytic structure around $(Z,(x_j))$ in $T$ is determined by the isomorphism types of the $Z_i$. Furthermore, a sufficiently small neighbourhood of $(Z,(x_j))$ in $T$ is isomorphic to a product $\prod_{i=1}^k U_i$, where the analytic structure of $U_i$ is determined by the isomorphism type of $Z_i$.

We then have
\[
 \nu_T((Z,(x_j))) = \nu_{\prod U_i}((Z, (x_j))) = \prod_i \nu_{U_i}(\pi_i((Z,(x_j))).
\]
Here $\pi_i$ is the projection to $U_i$, and the factor $\nu_{U_i}(\pi_i(Z, (x_j))$ depends only on the isomorphism type of $Z_i$. This implies that the level sets of $\nu_T$ intersected with $T((Q_i);A)$ are geometric. It follows that the complete level sets of $\nu_T$ are geometric.

As we have $\csm(\nu_T) = c_M(T)$ by definition, we get
\[
 \csm(T) = \csm(1_T) = \csm(\nu_T) + \csm(1_T - \nu_T) = c_M(T) + \sum i \csm(T_i)
\]
where the sum is finite and the $T_i$ are geometric subsets of $X^{[\![n]\!]}$ of lower dimension than $T$, with type depending only on the type of $T$. By induction on $\dim T$, we may assume the terms $\csm(T_i)$ are universal, and the hypothesis of the lemma is that $c_M(T)$ is as well.
Hence $\deg F(E^{[\![n]\!]}, \tsn) \cap \csm (T)$ is universal.
The functorial property of CSM classes implies
\[
 p_*(\csm(T)) = \deg(T/P) \csm(P) + \sum i \csm(P_i),
\]
where the $P_i$ are subsets of $X^{[n]}$, of lower dimension than $P$. It is not hard to show that the $P_i$ are geometric of type depending only on the type of $P$. Induction on the dimension of $P$ now gives the main theorem as claimed.
\end{proof}

\subsection{Partitions}
By a partition of $\{1, \ldots, n\}$, we mean a set $\alpha$ of disjoint subsets of $\{1, \ldots, n\}$, such that $\bigcup_{A \in \alpha} A = \{1, \ldots, n\}$.
Following \cite{li}, we will define schemes $X^{[\![\alpha]\!]}$ approximating $X^{[\![n]\!]}$ for each such $\alpha$.
We fix some notation and conventions with respect to partitions.

\begin{ndefn}
We let $\sim_\alpha$ be the equivalence relation on $\{1, \ldots, n\}$ given by letting the elements of $\alpha$ form equivalence classes.

We define a partial ordering on the set of partitions of $\{1, \ldots, n\}$ by letting $\alpha \le \beta$ if every element of $\alpha$ is contained in an element of $\beta$. Equivalently, $\alpha \le \beta$ if $\sim_\alpha$ is a finer relation than $\sim_\beta$.

We denote by $\Lambda$ the maximal partition under this ordering, that is, $\Lambda = \{\{1, \ldots ,n\}\}$. Given two partitions $\alpha, \beta$, we denote by $[\alpha, \beta]$ the set of partitions $\gamma$ such that $\alpha \le \gamma \le \beta$, and define $[\alpha,\beta)$ etc.\ similarly.
\end{ndefn}

\subsection{Approximating constructions}
\label{subsection:approximatingspaces}
From this point on we fix a closed, irreducible, geometric subscheme $T\subseteq X^{[\![n]\!]}$. In this section, we define the schemes $X^{[\![\alpha]\!]}$, the bundles $E^{[\![\alpha]\!]}$ and $T_{X^{[n]}}^{[\![\alpha]\!]}$, the subsets $T_\alpha \subseteq X^{[\![\alpha]\!]}$, and the cohomology classes $C_\alpha$ and $D$.

\begin{ndefn}
If $\alpha$ is a partition of $\{1, \ldots, n\}$, define the scheme $X^{[\![\alpha]\!]}$ by
\[
 X^{[\![\alpha]\!]} = \prod_{A\in\alpha}X^{[\![A]\!]},
\]
where $X^{[\![A]\!]} \cong X^{[\![|A|]\!]}$ and parametrises pairs $(Z, (p_i)_{i \in A})$ such that $\sum_{i \in A} p_i$ is the fundamental cycle of $Z$. 
\end{ndefn}
There is a natural morphism $X^{[\![\alpha]\!]} \to X^n$ defined by the decomposition $X^n = \prod_{A \in \alpha} X^A$ and the natural morphisms $X^{[\![A]\!]} \to X^A$. 

\begin{ndefn}
Define the vector bundles $E^{[\![\alpha]\!]}$ and $T_{X^{[n]}}^{[\![\alpha]\!]}$ on $X^{[\![\alpha]\!]}$ by
\[
 E^{[\![\alpha]\!]} = \bigoplus_{A\in\alpha} E^{[\![A]\!]},\ \ \ \ T_{X^{[n]}}^{[\![\alpha]\!]} = \bigoplus_{A\in\alpha} T_{X^{[A]}}^{[\![A]\!]},
\]
where we suppress pullback along the projection $X^{[\![\alpha]\!]} \to X^{[\![A]\!]}$.
\end{ndefn}

\begin{ndefn}
\label{def:diagonal}
If $\alpha$ is a partition of $\{1, \ldots, n\}$, denote by $\Delta_\alpha$ the subset of $X^n$ given by
\[
 \Delta_\alpha = \{\left(x_1, \ldots, x_n\right) \in X^n\mid x_i = x_j \text{ if } i\sim_\alpha j\}.
\]
We refer to the sets $\Delta_\alpha$ as \emph{diagonals}.
\end{ndefn}
Write $T = \ol{T((Q_i); A)}$, as is possible by Lemma \ref{lemma:geometric} (vi). 
We define the partition $\mu$ of $\{1, \ldots, n\}$ by $\mu = \{A_1, \ldots, A_k\}$.
The image of $T$ under $X^{\dbl n \dbr}$ is $\Delta_\mu$.

Our next task is to define schemes $T_\alpha \subseteq X^{[\![\alpha]\!]}$ birational to $T$, for those $\alpha$ such that this is possible. For any $\alpha$, let $f_\alpha : X^{[\![n]\!]} \dashrightarrow X^{[\![\alpha]\!]}$ be the natural isomorphism defined on the open set where the moduli problems $X^{[\![n]\!]}$ and $X^{[\![\alpha]\!]}$ solve are the same. Specifically, $f_\alpha$ is defined on the set of points $(Z, (x_i))$ where $x_i \not= x_j$ if $i \not\sim_\alpha j$. 

The locus where $f_\alpha$ is defined intersects $T$ if and only if $\alpha \ge \mu$. For such $\alpha$ we let $T_\alpha$ be the closure of $f_\alpha(T)$ in $X^{[\![\alpha]\!]}$.

Recall that $\wt{T_\alpha} \to T_\alpha$ denotes the Nash blow up. As $\wt{T_\alpha} \to T_\alpha$ is birational, the map $f_\alpha$ induces a natural rational map $g_\alpha : T \dashrightarrow \wt{T_\alpha}$. 
\begin{ndefn}
\label{defQ}
Let
\[
g := \left(g_\alpha\right)_{\alpha \ge \mu}: T \dashrightarrow \prod_{\alpha \ge \mu} \wt{T_\alpha},
\]
and define $Q$ to be the closure of $g(T)$ in $\prod \wt{T_\alpha}$.
\end{ndefn}

For every $\alpha$ there are birational proper morphisms $Q \to \wt{T_\alpha} \to T_\alpha$. Any cohomology class on $T_\alpha$ and $\wt{T_\alpha}$ may be pulled back along these morphisms without changing the degree, and we will suppress such pullbacks in the notation.

\subsection{Approximations of the cohomology classes}
\label{subsection:approximatingclasses}
The schemes and bundles defined in the previous section give rise to cohomology classes approximate to the one we want to compute (cf.\ \cite[Sec.\ 5.4]{li}). 
Recall that $\Lambda$ denotes the maximal partition of $\{1, \ldots, n\}$, and that $T_{\wt{T_\alpha}}$ is the Nash bundle on $\wt{T_\alpha}$.
\begin{ndefn}
\label{def:calpha}
Let $\alpha$ be a partition $\ge \mu$. Define the class $C_\alpha \in H^*(Q)$ by
\[
 C_\alpha = F\left(E^{[\![\alpha]\!]}, T_{X^{[n]}}^{[\![\alpha]\!]}\right) \cup c_\bullet\left(T_{\wt{T_\alpha}}\right).
\]
We let $C = C_\Lambda$.
\end{ndefn}
Note that the main theorem is reduced to the claim that $\deg (C \cap [Q])$ is universal.

\begin{ndefn}
\label{def:d}
Let $\alpha$ be a partition $\ge \mu$. Define the class $D_\alpha \in H^*(Q)$ by putting $D_\mu = C_\mu$, and for $\alpha > \mu$ let $D_\alpha$ be defined inductively by
\[
D_\alpha = C_\alpha - \sum_{\gamma \in \left[\mu, \alpha\right)} D_\gamma.
\]
We let $D = D_\Lambda$.
\end{ndefn}
\begin{remark}
It is not hard to show that in fact $D = \sum_{\alpha \ge \mu} (-1)^{|\alpha| -1} (|\alpha|-1)! C_\alpha$. Except in the proof of Proposition \ref{genFunProp}, we will not need this, and we work instead directly with the inductive definition of $D$.
\end{remark}

The motivation behind the definition of $D_\alpha$ is as follows. Firstly, it follows directly from the definition that if $\deg(D\cap[Q])$ and $\deg(C_\alpha\cap[Q])$ are universal for $\alpha \not = \Lambda$, then $\deg(C\cap[Q])$ is universal as well. 
Using induction on $n$, we will show in Lemma \ref{product} that the the degree of $C_\alpha\cap[Q]$ is universal for $\alpha \not= \Lambda$, reducing the problem to that of computing $D$. Secondly, $D$ is such that the restriction of $D$ to $Q \sm (Q|_{\Delta_\Lambda})$ vanishes, which allows us to reduce the computation of its degree to studying a small neighbourhood of $Q|_{\Delta_\Lambda}$.

\subsection{Reduction to $\deg(D \cap[Q])$}
If $\alpha$ is a nonmaximal partition of $\{1, \ldots, n\}$, the scheme $X^{[\![\alpha]\!]}$ is by definition a product of schemes $X^{[\![m]\!]}$ with $m < n$. Hence we can reduce the computation of $C_\alpha$ to computations of cohomology classes on such $X^{[\![m]\!]}$, and if these are universal, then $C_\alpha$ will be as well. This argument leads to the following induction result.
\begin{nlemma}
\label{product}
Let $m$ be a positive integer. Suppose that Theorem \ref{mainthm}/\ref{surfacethm} holds for every $n < m$, and suppose that for $n = m$ the degree of $D \cap [Q]$ is given by a universal linear polynomial in the Chern numbers of $(X,E)$. Then Theorem \ref{mainthm}/\ref{surfacethm} holds for $n = m$.
\end{nlemma}
\begin{proof}
Assume that the theorem holds for every $n < m$. We shall then show that for every partition $\alpha \in [\mu, \Lambda)$, the degree of $C_\alpha \cap [Q]$ is expressed by a universal polynomial. Since we have $C = D_\Lambda + \sum_{\alpha \in [\mu,\Lambda)} k_\alpha C_\alpha$, the statement of the lemma follows.

Let $\alpha \in [\mu, \Lambda)$. Recall first that by definition there is a product decomposition
\[
 X^{[\![\alpha]\!]} = \prod_{A\in \alpha} X^{[\![A]\!]}.
\]
This gives rise to a product decomposition $T_\alpha = \prod_{A\in \alpha} T_A$, where the $T_A \subseteq X^{[\![A]\!]}$ are closed, irreducible, geometric subsets. Since the Nash blow up preserves products, we have $\wt{T_\alpha} = \prod \wt{T_A}$, as well as bundle decompositions $E^{[\![\alpha]\!]} = \oplus E^{[\![A]\!]}$, $T_{\wt{T_\alpha}} = \oplus T_{\wt{T_A}}$ and $T_{X^{[n]}}^{[\![\alpha]\!]} = \oplus T_{X^{[n]}}^{[\![A]\!]}$.

Now, using the Whitney sum formula we can find an expression for
\[
 C_\alpha = F\left(E^{[\![\alpha]\!]}, T_{X^{[n]}}^{[\![\alpha]\!]}\right) \cdot c_\bullet\left(T_{\wt{T_\alpha}}\right)
\]
as a polynomial in the Chern classes of $E^{[\![A]\!]}$, $T_{X^{[n]}}^{[\![A]\!]}$, and $T_{\wt{T_A}}$ for different $A \in \alpha$. Since $\alpha < \Lambda$, we have $|A| < m$ for every $A \in \alpha$. By the induction hypothesis, we thus get a universal polynomial for
\[
\deg C_\alpha \cap [Q] = \deg C_\alpha \cap [T_\alpha],
\]
as required.

The claim about the degree of the universal polynomial $G$ in the main theorem also follows by induction, using the assumption that $\deg(D\cap[Q])$ is linear as a polynomial in the Chern numbers of $(X,E)$.
\end{proof}

Since the theorem is clear for $n = 0$, it now suffices to show that the degree of $D \cap [Q]$ is given by a linear polynomial in the Chern numbers of $(X,E)$.

\subsection{Relative constructions}
\label{relativesection}
We will show in Lemma \ref{rel} that the class $D$ vanishes when restricted to the part of $Q$ lying over the complement of the small diagonal $\Delta_\Lambda \subset X^n$. It may thus essentially be computed by looking at a neighbourhood of $Q|_{\Delta_\Lambda}$. The next step is now to use this to show that the degree of $D$ equals that of a class $\cD \in H^*(\cQ)$, where $\cQ$ is a scheme defined similarly to $Q$, but with $X^{[n]}$ replaced with the relative Hilbert scheme $\Hilb^n(\ol{TX}/X)$.

We therefore repeat the constructions of approximating schemes and classes in this relative setting. 
These are for the most part straightforward adaptations of the constructions in Sections \ref{subsection:approximatingspaces} and \ref{subsection:approximatingclasses}. 
The exception to this is the scheme $\cT$ that corresponds to $T$ (and so the schemes $\cT_\alpha$ and $ \cQ$ which are derived from $\cT$), where we impose the condition that the first marked point must lie in the 0-section $X \subset \ol{TX}$.

In order to integrate cohomology classes it will be convenient to work with proper schemes. Hence we let $\overline{TX}$ denote the $\PP^d$-bundle $\PP{\left(\cO_X \oplus TX\right)}$, with the convention that $\PP{\left(V\right)}$ is the set of lines through the origin in $V$. Let $\pi: \ol{TX} \to X$ be the projection, and let $\ol{TX}^n$, $\Sym^n(\ol{TX}/X)$ and $\Hilb^n(\ol{TX}/X)$ denote respectively the fibre product, relative symmetric product and relative Hilbert scheme of $\ol{TX}$ over $X$.

\begin{ndefn}
Define the scheme $\ol{TX}^{[\![n]\!]}$ by the cartesian diagram
\[
 \begin{CD}
  \ol{TX}^{[\![n]\!]} @>>> \Hilb^n\left(\ol{TX}/X\right) \\
  @VVV                   @VVV             \\
  \ol{TX}^n           @>>> \Sym^n\left(\ol{TX}/X\right).
 \end{CD}
\]
\end{ndefn}

\begin{ndefn}
Let $\cE^{[n]}$ be the tautological bundle on $\Hilb^n\left(\ol{TX}/X\right)$ corresponding to the vector bundle $\pi^*(E)$ on $\ol{TX}$, and let $\cE^{[\![n]\!]}$ be the pullback of $\cE^{[n]}$ to $\ol{TX}^{[\![n]\!]}$. If $\dim X = 1$ or 2, let $T_{\ol{TX}^{[n]}}$ be the relative tangent bundle of $\Hilb^{n}(\ol{TX}/X) \to X$, and denote its pullback to $\ol{TX}^{[\![n]\!]}$ by $T_{{\ol{TX}}^{[n]}}^{[\![n]\!]}$.
\end{ndefn}
\begin{ndefn}
Let $\alpha$ be a partition of $\{1, \ldots, n\}$. Define the scheme $\overline{TX}^{[\![\alpha]\!]}$ by
\[
 \ol{TX}^{[\![\alpha]\!]} = \prod_{A\in \alpha} \ol{TX}^{[\![A]\!]},
\]
where the product is the fibre product over $X$. Define the bundles $\cE^{[\![\alpha]\!]}$ and $\relt^{[\![\alpha]\!]}$ on $\ol{TX}^{[\![\alpha]\!]}$ by
\[
 \cE^{[\![\alpha]\!]} = \bigoplus_{A\in \alpha} \cE^{[\![A]\!]},\ \ \ \ \relt^{[\![\alpha]\!]} = \bigoplus_{A \in \alpha} T_{\ol{TX}^{[A]}}^{[\![A]\!]}
\]
suppressing notation for the natural pullbacks.
\end{ndefn}

Let $\cT^\pr$ be the subset of $\ol{TX}^{[\![n]\!]}$ consisting of the labelled subschemes $(Z, (x_i))$ such that $(Z,(x_i))$ is isomorphic to a labelled subscheme in $T$, where isomorphism of labelled schemes is as defined above Lemma \ref{lemma:geometric}. Equivalently, if $\ol{TX_x}^{[\![n]\!]}$ denotes the fibre of $\ol{TX}^{[\![n]\!]} \to X$ at $x$, then $\cT^\pr$ is the set such that $\cT^\pr \cap \ol{TX_x}^{[\![n]\!]}$ is a geometric subset of the same type as $T$, for every $x \in X$.

Let $r: \ol{TX}^{[\![n]\!]} \to \ol{TX}$ be the morphism defined by $r((Z, (x_i)) = x_1$, and let $\cT = r^{-1}(X) \cap \cT^\pr$, where $X \subset \ol{TX}$ is embedded by the 0-section.

For every partition $\alpha$, there is a rational map $f_\alpha : \ol{TX}^{[\![n]\!]} \dashrightarrow \ol{TX}^{[\![\alpha]\!]}$ defined where the moduli problems the two schemes solve are the same. Using these maps, we may replace $T$ by $\cT$ in Definition \ref{defQ} and the preceding paragraphs, thus defining the schemes $\cT_\alpha$ (for $\alpha \ge \mu$), $\wt{\cT_\alpha}$ and $\cQ$. We omit the details.

Finally, we define the relative analogues of the classes $C_\alpha$, $D_\alpha$.
\begin{ndefn}
Let $\alpha$ be a partition $\ge \mu$. Define the class $\cC_\alpha \in H^*(\cQ)$ by
\[
\cC_\alpha = F\left(\cE^{[\![\alpha]\!]}, \relt^{[\![\alpha]\!]}\right) \cup c_\bullet\left(T_{\wt{\cT_\alpha}}\right).
\]
We let $\cC = \cC_\Lambda$.
\end{ndefn}

\begin{ndefn}
Let $\alpha$ be a partition $\ge \mu$. Define the class $\cD_\alpha \in H^*(\cQ)$ by putting $\cD_\mu = \cC_\mu$, and for $\alpha > \mu$ let $\cD_\alpha$ be defined inductively by
\[
\cD_\alpha = \cC_\alpha - \sum_{\gamma \in [\mu,\alpha)} \cD_\gamma.
\]
We let $\cD = \cD_\Lambda$.
\end{ndefn}

\subsection{Relating $D_\Lambda$ to $\cD_\Lambda$}
\label{sec:relatingClasses}
Let $Q_0 \subseteq Q$ be the restriction of $Q$ to the small diagonal $X = \Delta_\Lambda \subset X^n$. Similarly let $\cQ_0\subseteq \cQ$ be the restriction of $\cQ$ to the set $X \subset \ol{TX}^n$ where the inclusion of $X$ is given by $n$ copies of the 0-section. The classes $D$ and $\cD$ are related by the following lemma and its corollary.
\begin{nlemma}
There exists a pair of open neighbourhoods $U^\prime \subset U$ of $Q_0$ in $Q$, a pair of open neighbourhoods $\cU^\prime \subset \cU$ of $\cQ_0$ in $\cQ$, a homeomorphism $\left(U^\prime, U\right) \to \left(\cU^\prime, \cU\right)$, and a class $\overline{\cD} \in H^*(\cQ, \cQ \setminus \cU^\prime)$ lifting $\cD\in H^*(\cQ)$, such that the composition
\[
 H^*(\cQ, \cQ \setminus \cU^\prime) \to H^*(\cU, \cU \setminus \cU^\prime) \to H^*(U, U \setminus U^\prime) \to H^*(Q, Q \setminus U^\prime) \to H^*(Q)
\]
sends the class $\overline{\cD}$ to $D$.
\label{rel}
\end{nlemma}
\begin{ncor}
 The degree of $D \cap [Q]$ equals the degree of $\cD \cap [\cQ]$.
\label{equaldegrees}
\end{ncor}
\begin{proof}[Proof of corollary]
 There are relative fundamental classes $[(\cQ,\cQ\sm\cU^\pr)]$, $[(\cU, \cU\sm \cU^\pr)]$, $[(U,U\sm U^\pr)]$ and $[(Q,Q\sm U^\pr)]$ in the appropriate homology groups. Replacing $H^*$ with $H_*$ in the above sequence (and reversing the arrows) each fundamental class is sent to the next, so in the composition the class $[Q]$ is sent to $[(\cQ,\cQ \sm \cU)]$. This implies
\[
 \deg \left(D \cap [Q]\right) = \deg \left(\ol{\cD} \cap [(\cQ,\cQ \sm \cU^\pr)]\right).
\]
Now, $[(\cQ, \cQ \sm \cU^\pr)]$ is the push-forward of $[\cQ]$, which shows that $\deg \left(\ol{\cD} \cap [(\cQ,\cQ \sm \cU^\pr)]\right) = \deg \left(\cD \cap [\cQ]\right)$, completing the proof.
\end{proof}

The proof of Lemma \ref{rel} is quite technical and is postponed to Section \ref{section:lemmaproof}. We now show how the main theorem follows from Corollary \ref{equaldegrees}.
\begin{proof}[Proof of main theorem]
By Corollary \ref{equaldegrees}, if $\deg \cD \cap [\cQ]$ is given by a universal linear polynomial, the same is true for $\deg D \cap [Q]$, which by Lemma \ref{product} would imply the main theorem.

Every construction made in Section \ref{relativesection} starting from $\ol{TX} \to X$ can be carried out with the bundle $TX \to X$ replaced by an arbitrary algebraic rank $d$ vector bundle. In particular, we may construct the analogue of $\cQ$ starting from the universal rank $d$ subbundle $U \to \Gr\left(d,N\right)$, where $N$ is any integer large enough that $TX$ is the pullback of $U$ along a continuous classifying map $X \to \Gr(d,N)$. Call this scheme $\cQ_{\Gr}$, let $q: \cQ_{\Gr} \to \Gr(d,N)$ be the natural morphism, and denote the analogues of $\cT_\alpha$ by $\cT_{\alpha,\Gr}$. 

Let $r$ be the rank of $E$, and let $M$ be a sufficiently large integer. There is then a Cartesian diagram
\[
 \begin{CD}
  \cQ @>f>> \cQ_{\Gr} \times \Gr(r, M)  \\
  @VpVV     @Vq\times \id VV                    \\
  X   @>g>> \Gr\left(d,N\right) \times \Gr(r, M)
 \end{CD}
\]
in the category of topological spaces, where $g$ is the product of the topological classifying maps for the bundles $TX$ and $E$. Note that $f$ and $g$ are in general not analytic.

Let $V \to \Gr(r,M)$ be the universal subbundle. If $\alpha \ge \mu$, let the bundle $V^{[\![\alpha]\!]}$ on $\cQ_{\Gr}\times \Gr(r,M)$ be defined by
\[
 V^{[\![\alpha]\!]} = V \otimes \cO_{\ol{U}^n}^{[\![\alpha]\!]},
\]
where $V$ and $\cO_{\ol{U}^n}^{[\![\alpha]\!]}$ are pulled back from $\Gr(r,M)$ and $\cQ_{\Gr}$, respectively. We then have $f^*(V^{[\![\alpha]\!]}) = \cE^{[\![\alpha]\!]}$. The scheme $\cQ_{\Gr}$ also carries a bundle $T_{\ol{U}^{[n]}}^{[\![\alpha]\!]}$, defined in the same way as $\relt^{[\![\alpha]\!]}$, and we have $f^*(T_{\ol{U}^{[n]}}^{[\![\alpha]\!]}) = \relt^{[\![\alpha]\!]}$.

For any $\alpha\ge \mu$, define the relative Nash bundle $T_{\wt{\cT}_{\alpha,\Gr}/\Gr}$ on $\cQ_{\Gr}$ as the kernel of the natural map $T_{\wt{\cT}_{\alpha,\Gr}} \to q^*\left(T_{\Gr(d,N)}\right)$. This map is surjective since $\cT_{\alpha,\Gr} \to \Gr(d,N)$ is a locally trivial fibration. Hence there is a short exact sequence
\[
 0\to T_{\wt{\cT}_{\alpha,\Gr}/\Gr} \to T_{\wt{\cT}_{\alpha,\Gr}} \to q^*\left(T_{\Gr(d,N)}\right)\to 0.
\]
Similarly define $T_{\wt{\cT}_{\alpha}/X}$ by the short exact sequence
\[
 0 \to T_{\wt{\cT}_{\alpha}/X} \to T_{\wt{\cT_\alpha}} \to p^*\left(T_X\right) \to 0.
\]

We then have $f^*(T_{\wt{\cT}_{\alpha,\Gr}/\Gr}) = T_{\wt{\cT}_{\alpha}/X}$. Define a bundle $G_\alpha$ on $\cQ_{\Gr}$ by $G_\alpha = T_{\wt{\cT}_{\alpha,\Gr}/\Gr} \oplus q^*\left(U\right)$. In K-classes we then have $f^*\left(G_\alpha\right) = T_{\wt{\cT_\alpha}}$. Define the class $\cC_{\alpha,\Gr} \in H^*(\cQ_{\Gr} \times \Gr(r,M))$ by
\[
 \cC_{\alpha,\Gr} = F\left(\cE_{\Gr}^{[\![\alpha]\!]}, T_{\ol{U}^{[n]}}^{[\![\alpha]\!]}\right) \cdot c_\bullet(G_\alpha).
\]
The above discussion shows that $f^*\left(\cC_{\alpha, \Gr}\right) = \cC_{\alpha}$. 

There are Gysin maps $p_!$ and $(q\times \id)_!$ in cohomology, defined by 
\[
p_!(\alpha) = PD(p_*(\alpha \cap [\cQ]))
\]
and 
\[
(q\times \id)_!(\alpha) = PD((q\times \id)_*(\alpha\cap [\cQ_{\Gr} \times \Gr(r,M)])),
\]
where PD denotes the Poincaré dual. These satisfy the relation $p_!f^* = g^*(q\times \id)_!$.

As a consequence, we get $p_!(\cC_{\alpha}) = g^*((q\times \id)_!(\cC_{\alpha,\Gr}))$, which implies
\begin{equation}
\label{eqn:finalEquation}
 \deg(\cC_\alpha \cap [\cQ]) = \deg(g^*((q\times \id)_!(\cC_{\alpha,\Gr})) \cap [X]).
\end{equation}

Any rational cohomology class on $\Gr\left(d,N\right) \times \Gr(r,M)$ can be expressed as a polynomial in Chern classes of the universal bundles. 
If moreover $M$ and $N$ are sufficiently large with respect to $d$, there are no relations between these Chern classes in degree $2d$, so this polynomial is unique. In particular, the degree $2d$ part of $(q\times\id)_!\left(\cC_{\alpha,\Gr}\right)$ is equal to such a polynomial. This polynomial is independent of $M$ and $N$, because the class $\cC_{\alpha,\Gr}$ is preserved by the pullbacks induced by the natural morphisms $\Gr(d,N) \to \Gr(d,N+1)$ and $\Gr(r,M) \to \Gr(r,M+1)$.

It follows then that the right hand side of \ref{eqn:finalEquation} is equal to a linear combination of the Chern numbers of $(X,E)$. Consequently, $\deg \cC_\alpha \cap [\cQ]$ is a universal linear combination of the Chern numbers of $(X,E)$. As $\cD$ is a linear combination of the $\cC_\alpha$, the same is true for $\deg \cD \cap [\cQ]$, which is what we needed to show.
\end{proof}

\section{Proof of Lemma \ref{rel}}
\label{section:lemmaproof}
\subsection{Defining the map from $Q$ to $\cQ$}
\label{defMap}
Let $p_1, p_2 : X\times X \to X$ be the projection to the first and second factor, and let $\pi: TX \to X$ be the tangent bundle.
\begin{nlemma}
\label{f1}
 There is an open neighbourhood $U_1$ of the diagonal $\Delta \subset X\times X$, an open neighbourhood $\cU_1$ of the 0-section $X\subset TX$ and a homeomorphism $f_1:U_1 \to \cU_1$, such that
\[
 \pi \circ f_1 = p_1
\]
and such that $f_1|_\Delta$ is the identitification between $\Delta$ and the 0-section of $TX$. Furthermore, the restriction of $f_1$ to each fibre $p_1^{-1}(x)$ is holomorphic. 

There is an isomorphism of topological vector bundles $p_1^*(E)|_U \to p_2^*(E)|_U$, which is an isomorphism of holomorphic bundles on the restriction to each fibre $p_1^{-1}(x)$.
\label{Li}
\end{nlemma}
\begin{proof}
See \cite[Lemma 2.4]{li} for the first two statements. Holomorphic exponential maps can be constructed on small open sets, and these can be globalised using a partition of unity. 
This globalising preserves holomorphicity on fibres of $p_1$ as required.

For the statement about $E$, we argue similarly. Cover $X$ with open balls $B_i$, and choose holomorphic trivialisations $E|_{B_i} \cong \cO_{B_i}^n$. Using these, define local isomorphisms $g_i: p_1^*(E)|_{B_{i} \times B_{i}} \to p_2^*(E)|_{B_{i} \times B_{i}}$. 
Let $\{t_i\}$ be a smooth partition of unity subordinate to the covering $\{B_i\}$, and define $g : p_1^*(E) \to p_2^*(E) $ at $x \in U_1$ by $g(x)  = \sum t_i(p_1(x)) \cdot g_i$. 

This map $g$ is holomorphic on fibres of $p_1$. Restricted to $\Delta$, the map $g$ is the identity, and so after shrinking $U_1$ if necessary, $g$ is an isomorphism.
\end{proof}

Let $X^{[\![n]\!]}_0 \subset X^{[\![n]\!]}$ be the set of pairs $(Z,(x_i))$ such that $Z$ is supported at a single point, and let $\overline{TX}^{[\![n]\!]}_0 \subset \ol{TX}^{[\![n]\!]}$ denote the set of pairs $(Z,(x_i))$ such that $Z$ is supported at the 0-section of $\ol{TX}$. Let $q: X^{[\![n]\!]} \to X$ be defined by $q(Z,(x_i)) = x_1$, and let $r : \ol{TX}^{[\![n]\!]} \to X$ be defined by $r((Z, (x_i))) = \pi(x_1)$. Let $W$ be the set of pairs $(Z, (x_i)) \in \ol{TX}^{[\![n]\!]}$ such that $x_1$ lies in the 0-section of $\ol{TX}$.
\begin{nlemma}
\label{f2}
There is an open neighbourhood $U_2$ of $X^{[\![n]\!]}_0$ in $X^{[\![n]\!]}$, an open neighbourhood $\cU_2$ of $\overline{TX}^{[\![n]\!]}_0$ in $W$, and a homeomorphism $f_2 : U_2 \to \cU_2$, such that $q = r \circ f_2$.
Furthermore, $f_{2}|_{q^{-1}(x)}$ is holomorphic for all $x \in X$. 
There are isomorphisms of topological vector bundles $f_2^*\left(T_{\ol{TX}^{[n]}}^{[\![n]\!]}\right) \to T_{X^{[n]}}^{[\![n]\!]}$ and $f_2^*\left(\cE^{[\![n]\!]}\right) \to E^{[\![n]\!]}$.
\end{nlemma}
\begin{proof}
Let $U_1$, $\cU_1$ and $f_1$ be as provided by Lemma \ref{f1}. 
Define $f_{2}$ by
\[
f_{2}((Z,(x_{i})) = ((f_{1})_{*}(\{q(x)\} \times Z), (f_{1}(q(x), x_{i})).
\]
The right hand side is well-defined if $\{q(x)\} \times Z$ is contained in $U_{1}$.
Let $U_{2} \subset X^{\dbl n \dbr}$ be the open set where this is the case, so that $f_{2}$ is defined on $U_{2}$.
Then $f_{2}$ is a local homeomorphism such that $q = r \circ f_{2}$, and $f_{2}|_{q^{-1}(x)}$ is analytic for all $x \in X$.

We now describe the isomorphism $T^{\dbl n\dbr}_{X^{[n]}} \cong f_2^*(T_{\ol{TX}^{[n]}}^{[\![n]\!]})$.
Over a point $(Z, (x_i)) \in X^{\dbl n \dbr}$, this is the composition
\begin{align*}
\left( T^{\dbl n \dbr}_{X^{[n]}} \right)_ {(Z, (x_i))} = \left ( T_{X^{[n]}}\right )_Z \stackrel{\cong}{\to} \left (T_{\ol{TX}^{[n]}}\right )_{(f_1)_*(\{x_1\} \times Z)} = \left(T^{\dbl n \dbr}_{\ol{TX}^{[n]}}\right )_{f_2(Z,(x_i))},
\end{align*}
where the middle map is the differential of the map $Z \mapsto (f_1)_*(\{x_1\} \times Z)$.

Finally, we describe the isomorphism $f_2^*(\cE^{\dbl n \dbr}) \cong E^{\dbl n \dbr}$.
Over a point $(Z, (x_i)) \in X^{\dbl n \dbr}$, we have
\begin{align*}
f_2^*\left(\cE^{\dbl n \dbr}\right)_{(Z, (x_i))} &= E_{x_1} \otimes H^0\left(\cO_{f_1(\{x_1\} \times Z)}\right) \cong E_{x_1}\otimes H^0\left(\cO_{\{x_1\}\times Z}\right) \\
&\cong H^0\left(\{x_1\} \times Z, p_1^*(E)|_{\{x_1\} \times Z}\right).
\end{align*}
The isomorphism $p_1^*(E) \to p_2^*(E)$ of Lemma \ref{f1} now gives an isomorphism
\[
H^0\left(\{x_1\} \times Z, p_1^*(E)|_{\{x_1\} \times Z}\right) \stackrel{\sim}{\to}H^0\left(\{x_1\} \times Z, p_2^*(E)|_{\{x_1\} \times Z}\right) = E^{\dbl n \dbr}_{(Z, (x_i))}.
\]
\end{proof}

Recall that $Q_0 \subset Q$ and $\cQ_0 \subset \cQ$ are the subsets of points having image in $\Delta \subset X^n$ and $X \subset \ol{TX}^n$ under the natural morphisms $Q \to X^n$ and $\cQ \to \ol{TX}^n$, respectively.

Let the relative Nash bundles $T_{\wt{T_\alpha}/X}$ and $T_{\wt{\cT_\alpha}/X}$ be the kernels of the surjective homomorphisms $T_{\wt{T_\alpha}} \to q^{-1}(T_X)$ and $T_{\wt{\cT_\alpha}} \to r^{-1}(T_X)$, respectively.

\begin{nlemma}
\label{f3}
 There is an open neighbourhood $U$ of $Q_0$ in $Q$, an open neighbourhood $\cU$ of $\cQ_0$ in $\cQ$, and a homeomorphism $f : U \to \cU$, as well as isomorphisms of topological vector bundles
\[
 f^*\left(\cE^{[\![\alpha]\!]}\right) \to E^{[\![\alpha]\!]},
\]
\[
 f^*\left(T_{\ol{TX}^{[n]}}^{[\![\alpha]\!]}\right) \to T_{X^{[n]}}^{[\![\alpha]\!]},
\]
and
\[
 f^*\left(T_{\wt{\cT_\alpha}/X}\right) \to T_{\wt{T_\alpha}/X}.
\]
\end{nlemma}
\begin{proof}
The map $f_2$ constructed in Lemma \ref{f2} gives rise to local isomorphisms $T_\alpha \to \cT_\alpha$. 
The Nash blow up is determined analytically locally, and $T_\alpha \stackrel{q}{\to} X$ and $\cT_\alpha \stackrel{r}{\to} X$ are both locally trivial fibrations in a neighbourhood of $Q_0$ and $\cQ_0$.

If $F, U$ are complex analytic spaces with $U$ smooth, then the Nash blow-up $\wt{F \times U}$ is canonically isomorphic to $\wt{F} \times U$.
It follows that the Nash blow ups $\wt{T}_\alpha$ and $\wt{\cT}_\alpha$ are locally trivial fibrations over $X$, and that the local isomorphisms $T_\alpha \to \cT_\alpha$ extend uniquely to isomorphisms of the Nash blow-ups.
This in turn induces a local isomorphism $Q \to \cQ$. 

The first two bundle isomorphisms are induced by the ones produced in Lemma \ref{f2}, and the third follows similarly, taking into account the fact that $T_\alpha \to \cT_\alpha$ is holomorphic on the fibres of $q$ and $r$.
\end{proof}

\subsection{Construction of $\ol{D}$, $\ol{\cD}$}
In order to prove the claim of Lemma \ref{rel}, we will construct the relative cohomology classes $\ol{D} \in H^*(Q, Q \sm U)$ and $\ol{\cD} \in H^*(\cQ, \cQ \sm V)$ explicitly as singular cochains. Adapting the argument in Section 5.3 of \cite{li}, we first define certain open subsets $U_\alpha$, $V_{\alpha,\beta}$ of $X^n$ which we will later use to compare the bundles $E^{[\![\alpha]\!]}$ (or $T_{X^{[n]}}^{[\![\alpha]\!]}$, $T_{\wt{T}_\alpha}$) for various $\alpha$.

Let $d_X$ be a metric on $X$ inducing the Euclidean topology. Define the metric $d_{X^n}$ on $X^n$ by
\[
 d_{X^n}\left(\left(x_i\right)_{i=1}^n, \left(y_i\right)_{i=1}^n\right) = \max_{1\le i \le n} d_X\left(x_i, y_i\right).
\]
Let $d_{\ol{TX}}$ be a metric on $\ol{TX}$ inducing the Euclidean topology, and let $d_{\ol{TX}^n}$ be the metric on $\ol{TX}^n$ defined in the same way as $d_{X^n}$.

Let $x \in X^n$ and let $\alpha$ be a partition. Recall that $\Delta_\alpha \subseteq X^n$ is the diagonal set
\[
 \{(x_1, \ldots, x_n) \mid x_i = x_j\text{ if }i \sim_\alpha j\},
\]
and let $\Delta_\alpha^\pr \subset \ol{TX}^n$ be defined in the same way. In the following, we will use the inequalities
\begin{equation}
\label{keyinequality}
 \frac{1}{2}\sup_{i,j \mid i\sim_\alpha j} d_X\left(x_i,x_j\right) \le d_{X^n}\left(x, \Delta_\alpha\right) \le \sup_{i,j \mid i\sim_\alpha j} d_X\left(x_i, x_j\right)
\end{equation}
and their variants for $d_{\ol{TX}}$ and $d_{\ol{TX}^n}$, all of which follow easily from the definitions and the triangle inequality.

\begin{ndefn}
Let $P\left(n\right)$ be the set of partitions of $\{1, \ldots, n\}$, and let $\epsilon : P\left(n\right) \to \RR^{>0}$ be a function. At various points in the proof, the quantities
\[
 \max_{\alpha \in P\left(n\right)} \epsilon\left(\alpha\right)
\]
and
\[
 \max_{\alpha < \beta \in P\left(n\right)} \frac{\epsilon\left(\alpha\right)}{\epsilon\left(\beta\right)}
\]
will be assumed to be sufficiently small.
\end{ndefn}

\begin{ndefn}
Let $U_\alpha \subset X^n$ be the open $\epsilon\left(\alpha\right)$-neighbourhood of $\Delta_\alpha\subset X^n$, and let $\cU_\alpha$ be the open $\epsilon\left(\alpha\right)$-neighbourhood of $\Delta_\alpha^\pr \subset \ol{TX}^n$.
\end{ndefn}
\begin{ndefn}
Let $\alpha, \beta$ be partitions such that $\alpha < \beta$. Define the set $V_{\alpha,\beta} \subset X^n$ as
\[
 V_{\alpha,\beta} = \left(X^n\sm U_\beta\right) \sm \left(\bigcup_{\stackrel{\gamma < \beta}{\gamma \not\le \alpha}} \ol{U}_\gamma\right),
\]
and define the set $\cV_{\alpha,\beta} \subset \ol{TX}^n$ as
\[
 \cV_{\alpha,\beta} = \left(\ol{TX}^n\sm \cU_\beta\right) \sm \left(\bigcup_{\stackrel{\gamma < \beta}{\gamma \not\le \alpha}} \ol{\cU}_\gamma\right).
\]
\end{ndefn}

Let $i,j \in \{1, \ldots ,n\}$. Define an equivalence relation $\sim_{\left(i,j\right)}$ on the set of partitions by saying $\alpha \sim_{\left(i,j\right)}\beta$ if $\sim_\alpha$ and $\sim_\beta$ agree when evaluated on the pair $\left(i,j\right)$. For any pair of partitions $\alpha, \beta$, define the set $\Delta_{\alpha\beta} \subseteq X^n$ to be the set of points over which $X^{[\![\alpha]\!]}$ and $X^{[\![\beta]\!]}$ are not canonically equal. Explicitly, we have
\[
 \Delta_{\alpha\beta} = \bigcup_{i,j\mid \alpha \not \sim_{\left(i,j\right)} \beta} \Delta_{ij},
\]
where $\Delta_{ij}$ denotes the set of points $x \in X^n$ for which $x_i = x_j$. Define $\Delta_{\alpha\beta}^\pr \subset \ol{TX}^n$ similarly.

The following lemma summarises the important properties of $V_{\alpha,\beta}$ and $\cV_{\alpha,\beta}$.
\begin{nlemma}\mbox{}

\begin{enumerate}
\item[(i)] Let $\beta$ be a partition. The sets $\{V_{\alpha, \beta}\}_{\alpha< \beta}$ form an open covering of $X^n\setminus U_\beta$. The sets $\{\cV_{\alpha,\beta}\}_{\alpha < \beta}$ form an open covering of $\left(\overline{TX}^n\right)\setminus \cU_\beta$.
\item[(ii)] Let $\alpha, \beta,\gamma$ be partitions such that $\alpha < \beta$, $\gamma \le \beta$ and $\gamma \not \le \alpha$. Then 
\[
U_\gamma \cap V_{\alpha,\beta} = \emptyset\ \ \ \ \text{and}\ \ \ \ \cU_\gamma \cap \cV_{\alpha,\beta} = \emptyset.
\]
\item[(iii)] Let $\tau = \min_\gamma \epsilon(\gamma)$, and let $\alpha < \beta$. If $x \in V_{\alpha,\beta}$, we have $d(x, \Delta_{\alpha\beta}) \ge \tau$. If $x \in \cV_{\alpha,\beta}$, we have $d(x, \Delta^\pr_{\alpha\beta}) \ge \tau$.
\end{enumerate}
\label{opencovering}
\end{nlemma}
\begin{proof}
We prove the statements for $V_{\alpha,\beta}$; the case of $\cV_{\alpha,\beta}$ is exactly the same.

(i): Assume $x = (x_i) \in X^n \sm U_\beta$, and let $\alpha$ be maximal among partitions $< \beta$ such that $x \in \ol{U}_\alpha$. Such a partition exists because for the smallest partition $\omega$, we have $U_\omega = X^n$. We claim that $x\in V_{\alpha,\beta}$.

Assume $x\not \in V_{\alpha,\beta}$, there is then a partition $\gamma$ such that $\gamma \not \le \alpha$, $\gamma < \beta$ and $x \in \ol{U}_\gamma$. By the maximality property of $\alpha$, we cannot have $\alpha < \gamma$. It follows that $\alpha,\gamma < (\alpha\vee\gamma) \le \beta$, where $\alpha \vee \gamma$ is the smallest partition majorising $\alpha$ and $\gamma$.

Let $i,j$ be two indices such that $i\sim_{\alpha\vee\gamma} j$, and such that $d\left(x_i,x_j\right)$ is maximal for pairs with this property. There is a sequence of integers $i_1, i_2, \ldots, i_r$ such that $i_1 = i$, $i_r = j$ and such that for every $k$ with $1 \le k < r$, either $i_k \sim_\alpha i_{k+1}$ or $i_k \sim_\gamma i_{k+1}$ is true. By \eqref{keyinequality}, we now have
\[
 d\left(x, \Delta_{\alpha\vee\gamma}\right) \le d\left(x_i,x_j\right) \le d\left(x_{i_1}, x_{i_2}\right) + \cdots + d\left(x_{i_{r-1}}, x_{i_r}\right).
\]
Since $x \in \ol{U}_\alpha$, we have $d(x, \Delta_\alpha) \le \epsilon(\alpha)$, and similarly for $\gamma$. By \eqref{keyinequality}, each term in the above sum is therefore $\le 2\max(\epsilon(\alpha),\epsilon(\gamma))$. The sum is therefore smaller than
\[
2\left(r-1\right)\max(\epsilon(\alpha),\epsilon(\gamma)) < \epsilon\left(\alpha\vee\gamma\right),
\]
where the last inequality uses the second smallness assumption in the definition of $\epsilon$.
Hence we have $d\left(x, \Delta_{\alpha\vee\gamma}\right) < \epsilon\left(\alpha\vee\gamma\right)$, so that $x \in U_{\alpha\vee\gamma}$. If $\alpha \vee \gamma \not= \beta$, this contradicts the maximality of $\alpha$, and if $\alpha \vee \gamma = \beta$, it contradicts the assumption that $x \not \in U_\beta$.

(ii): Obvious from the definition.

(iii): Let $x = (x_i) \in V_{\alpha,\beta}$. For every $i, j \in \{1, \ldots, n\}$, let $\gamma_{i,j}$ be the partition defined by the equivalence relation such that $i \sim_{\gamma_{i,j}} j$ and no other nontrivial relations hold. If $\alpha < \beta$, we have
\[
 \Delta_{\alpha\beta} = \bigcup_{\substack{i \sim_\beta j \\ i \not\sim_\alpha j}} \Delta_{\gamma_{i,j}}.
\]
For every pair $i,j$ occurring in the union, we have $\gamma_{i,j} \not\le \alpha$, hence by part (ii) of the lemma we have $x \not\in U_{\gamma_{i,j}}$. This gives
\[
 d(x,\Delta_{\alpha\beta}) = \min_{\substack{i \sim_\beta j \\ i \not\sim_\alpha j}} d(x, \Delta_{\gamma_{i,j}}) > \min_{\substack{i \sim_\beta j \\ i \not\sim_\alpha j}} \epsilon(\gamma_{i,j}) \ge \tau.
\]
\end{proof}

Recall that $T_{\wt{T}_\alpha/X}$ and $T_{\wt{\cT}_\alpha/X}$ are the relative Nash bundles defined above Lemma \ref{f3}. Denote by $\cO^N$ the trivial bundle of rank $N$.
\begin{nlemma}\mbox{}
\begin{enumerate} 
\item[(i)] For each $\alpha \ge \mu$, let $F_\alpha$ denote either $E^{[\![\alpha]\!]}, T_{\wt{T_\alpha}/X}$ or $T_{X^{[n]}}^{[\![\alpha]\!]}$, considered as a topological bundle on $Q$. There is an integer $M$ and for every $\alpha$ an injection $i_\alpha : F_\alpha \to \cO^M$ such that if $\alpha \le \beta$, then over $V_{\alpha, \beta}$ the bundles $i_\alpha\left(E^{[\![\alpha]\!]}\right)$ and $i_\beta\left(E^{[\![\alpha]\!]}\right)$ are equal as subbundles of $\cO^M$.
 
\item[(ii)] For each $\alpha \ge \mu$, let $\cF_\alpha$ denote either $\cE^{[\![\alpha]\!]}, T_{\wt{\cT_\alpha}/X}$ or $\relt^{[\![\alpha]\!]}$, considered as a topological bundle on $\cQ$. There is an integer $N$ and for every $\alpha$ an injection $j_\alpha : \cF_\alpha \to \cO^N$ such that if $\alpha \le \beta$, then over $\cV_{\alpha, \beta}$ the bundles $j_\alpha\left(\cF^{[\![\alpha]\!]}\right)$ and $i_\beta\left(\cF^{[\![\alpha]\!]}\right)$ are equal as subbundles of $\cO^N$.

\item[(iii)] Let $f : U \to \cU$ be the local homeomorphism constructed in Lemma \ref{f3}. We may choose $M = N$ and the injections $i_\alpha, j_\alpha$ in such a way that the diagram
\[
\begin{CD}
 f^*\left(\cF_\alpha\right) @<\cong<< F_\alpha \\
 @Vf^*\left(j_\alpha\right)VV         @Vi_\alpha VV \\
 f^*\left(\cO^N\right)      @=        \cO^N
\end{CD}
\]
of bundles on $U$ commutes, where the upper isomorphism is the one constructed in Lemma \ref{f3}.
\end{enumerate}
\label{bundles}
\end{nlemma}
\begin{proof}
In this proof, all partitions are assumed to be $\ge \mu$.

(i): For each $\alpha$ choose first an injective homomorphism $i_\alpha^\pr : F_\alpha \to \cO^{M_\alpha}$. Let $k_{\alpha\beta} : F_\alpha \to F_\beta$ be the natural isomorphisms, defined over $X^n\sm \Delta_{\alpha\beta}$. Recall that $\tau = \min_\alpha \epsilon(\alpha)$. Let $t :\RR^{\ge 0} \to \RR^{\ge 0}$ be a continuous function such that when $x \ge \tau$ we have $t\left(x\right) = 1$ and such that when $x \le \tau/2$ we have $t\left(x\right) = 0$. For $x \in Q$, put $t_{ij}\left(x\right) = t\left(d\left(p(x),\Delta_{ij}\right)\right)$, where $p : Q \to X^n$ is the natural morphism. Let
\[
t_{\alpha\beta} = \prod_{i< j \mid \alpha \not\sim_{\left(i,j\right)} \beta} t_{ij}.
\]
As $t_{\alpha\beta}$ is supported away from $\Delta_{\alpha\beta}$, the homomorphism $t_{\alpha\beta}\cdot k_{\alpha\beta}$ is defined on the whole of $Q$. 

For any two partitions $\alpha$ and $\gamma$, let $i_{\alpha\gamma} = i_\gamma^\pr \circ \left(t_{\alpha\gamma} \cdot k_{\alpha\gamma}\right)$. We take $M = \sum M_\gamma$, and set
\[
 i_\alpha = \oplus_\gamma \left(i_{\alpha\gamma}\right)_\gamma : E_\alpha \to \oplus \CC^{M_\gamma} = \CC^M.
\]
It remains to show that $i_\alpha$ has the properties stated. As $i_{\alpha\alpha} = i_\alpha^\pr$ is injective, it is clear that $i_\alpha$ is an injection. 

Let $\alpha < \beta$, and let $\gamma$ be arbitrary. First we show that if $x \in Q$ lies over $V_{\alpha,\beta}$, then 
\begin{equation}
\label{claim1}
\left(t_{\alpha\gamma}\cdot k_{\alpha\gamma}\right)\left(x\right) = \left(t_{\beta\gamma} \cdot k_{\beta\gamma}\circ k_{\alpha\beta}\right)\left(x\right).
\end{equation}
To this end, observe that
\[
 \frac{t_{\alpha\gamma}}{t_{\beta\gamma}} = \prod_{i< j\mid \alpha \not\sim_{\left(i,j\right)} \beta} t_{ij}^{\delta(i,j)},
\]
where each $\delta(i,j) \in \{-1,0,1\}$. Now, since $p(x) \in V_{\alpha,\beta}$, we have $t_{ij}(x) = 1$ for every factor above, using Lemma \ref{opencovering} (iii). Hence, the equation $t_{\alpha\gamma}\left(x\right) = t_{\beta\gamma}\left(x\right)$ holds. If $t_{\alpha\gamma} = 0$, this shows \eqref{claim1}. If not, then all the morphisms $k_{\alpha\beta}, k_{\beta\gamma}, k_{\alpha\gamma}$ are defined at $x$, and by the naturality of these the cocycle condition $k_{\alpha\gamma} = k_{\beta\gamma}\circ k_{\alpha\beta}$ holds.

The above paragraph shows that $i_{\alpha\gamma} = i_{\beta\gamma} \circ k_{\alpha\beta}$ over $V_{\alpha,\beta}$ and hence $i_\alpha = i_\beta \circ k_{\alpha\beta}$. Consequently the two subbundles $i_\alpha\left(E_\alpha\right)$ and $i_\beta\left(E_\beta\right)$ of $\CC^M$ are equal as claimed.

\bigskip
(ii): Similar to (i).

\bigskip
(iii): Let $k_{\alpha\beta}^\pr : \cF_\alpha \to \cF_\beta$ be the homomorphisms defined like the $k_{\alpha\beta}$ in the proof of (i). Let $g_\alpha : F_\alpha \to f^*(\cF_\alpha)$ be the isomorphism of Lemma \ref{f3}. We then have $g_\beta \circ k_{\alpha\beta} = f^*\left(k_{\alpha\beta}^\pr\right)\circ g_\alpha$.

Let $U^\pr \subset U$ and $\cU^\pr \subset \cU$ be smaller open neighbourhoods of $Q_0$ and $\cQ_0$, such that $f(U^\pr) = \cU^\pr$. Let $s$ be a nonnegative real function which is 1 on $U^\pr$ and 0 on the complement of $U$. Replace $i_\alpha$ with
\[
 i_\alpha^{\pr} := \left(\left(1-s\right) i_\alpha, s\cdot f^*(j_\alpha)\circ g_\alpha\right) : F_\alpha \to \cO^{M + N},
\]
where $j_\alpha : \cF_\alpha \to \cO^N$ is the homomorphism of part \emph{(ii)}. Replacing $j_\alpha$ with $j_\alpha^\pr = (0,j_\alpha) \to \cO^M \oplus \cO^N$, then obviously $i_\alpha^{\pr} = f^*(j^\pr_\alpha) \circ g_\alpha$ over $U^\pr$. After replacing $U$ with $U^\pr$, it now remains to be shown that the statement in part (i) of the lemma still holds for $i_\alpha^{\pr}$.

After enlarging the metric on $\ol{TX}$ and shrinking $U$ and $\cU$, we may assume that if $x \in U$ lies over $V_{\alpha,\beta}$, then $f(x)$ lies over $\cV_{\alpha,\beta}$. Over $V_{\alpha,\beta}$ we thus have
\begin{align*}
 i_\alpha^{\pr} &= \left(\left(1-s\right) i_\alpha, s\cdot f^*(j_\alpha) \circ g_\alpha\right) = \left(\left(1-s\right) i_\beta \circ k_{\alpha\beta}, s\cdot f^*(j_\beta \circ k_{\alpha\beta}^\pr) \circ g_\alpha\right)\\
&= \left(\left(1-s\right) i_\beta, s\cdot f^*(j_\beta) \circ g_\beta\right) \circ k_{\alpha\beta} = i_\beta^\pr \circ k_{\alpha\beta}.
\end{align*}
\end{proof}

We are now in a position to finish the proof of Lemma \ref{rel}.
\begin{proof}[Proof of Lemma \ref{rel}:]
 Let $\alpha$ be a partition $\ge \mu$. The inclusions produced in Lemma \ref{opencovering} (ii) define a continuous map $f_\alpha : \cQ \to \Gr_1\times \Gr_2 \times \Gr_3$, where the $\Gr_i$ are Grassmannians with universal bundles $W_i$, such that
\[
 f_\alpha^*(W_1) = \cE^{[\![\alpha]\!]}\ \ \ \ f_\alpha^*(W_2) = \relt^{[\![\alpha]\!]}, \ \ \ \ f_\alpha^*(W_3) = T_{\wt{\cT_{\alpha}/X}} \oplus r^*(T_X).
\]
Note that at the level of Chern classes, we may interchange $f_\alpha^*(W_3)$ with $T_{\wt{\cT_{\alpha}}}$, as these are topologically isomorphic.

Choose a singular cochain $A$ on $\Gr_1 \times \Gr_2 \times \Gr_3$ representing the class
\[
 F(W_1, W_2) \cdot c_\bullet(W_3).
\]
Define singular cochains $\ol{\cC}_\alpha$ and $\ol{\cD}_\alpha$ by $\ol{\cC}_\alpha = f_{\alpha}^*(A)$ and
\[
 \ol{\cD}_\alpha = \ol{\cC}_\alpha - \sum_{\gamma \in [\mu, \alpha)} \ol{\cD}_\gamma.
\]
Clearly, the class of $\ol{\cC}_\alpha$ and $\ol{\cD}_\alpha$ is $\cC_\alpha$ and $\cD_\alpha$, respectively. Let $\ol{\cC} = \ol{\cC}_\Lambda$ and $\ol{\cD} = \ol{\cD}_\Lambda$.

We now claim that $\ol{\cD}|_{\ol{TX}^n\sm \cU_\Lambda} = 0$. To prove this, we show that $\ol{\cD}_\beta|_{\ol{TX}^n\sm \cU_\beta} = 0$ for any partition $\beta$ by ascending induction on the ordering of partitions. The base case is clear, as $\cQ|_{\ol{TX}^n\sm \cU_\mu} = \emptyset$. 

Assume now that $\ol{\cD}_\alpha|_{\ol{TX}^n\sm \cU_\alpha} = 0$ for every $\alpha < \beta$. We must show that for every singular $m$-simplex $a:\Delta^m \to \cQ|_{\ol{TX}^n \sm \cU_\beta}$ we have $\ol{\cD}_\beta(a) = 0$.

Since $\ol{\cD}_\beta$ is a cocycle, we may replace $a$ by any subdivision of $a$ and prove the vanishing for each simplex in the subdivision. By Lemma \ref{opencovering} (i), $\{\cV_{\alpha,\beta}\}_{\alpha < \beta}$ is an open covering of $\ol{TX}^n \sm \cU_\beta$, so we may assume there is an $\alpha < \beta$ such that $a$ maps into $\cQ|_{\cV_{\alpha,\beta}}$. 

If $\gamma < \beta$ is such that $\gamma \not \le \alpha$, then by Lemma \ref{opencovering} (ii) we have $\cU_\gamma \cap \cV_{\alpha,\beta}=\emptyset$, and so by the induction hypothesis $\ol{\cD}_\gamma\left(a\right) = 0$. This implies
\[
 \ol{\cD}_\beta\left(a\right) = \ol{\cC}_\beta\left(a\right) - \sum_{\gamma \le \alpha} \ol{\cD}_\alpha\left(a\right) = \ol{\cC}_\beta\left(a\right) - \ol{\cC}_\alpha\left(a\right),
\]
where the last equality follows directly from the definition of $\ol{\cD}_\alpha$. By Lemma \ref{bundles} (ii) we have $f_\alpha = f_\beta$ over $\cV_{\alpha,\beta}$, hence $\ol{\cD}_\beta\left(a\right) = 0$ and the claim follows.

Taking $\epsilon$ small enough we may assume that $\cQ\sm \cU^\pr$ lies over $\ol{TX}^n \sm \cU_\Lambda$, and then the above shows $\ol{\cD}$ is a relative cocycle for the pair $\left(\cQ, \cQ \sm \cU^\pr\right)$. The exact same construction performed on $Q$ produces a cochain $\ol{D}_{\alpha}$ of class $D_\alpha$ which vanishes on $Q \sm U^\pr$. By Lemma \ref{bundles} (iii), the homomorphism
\[
 H^*(\cQ,\cQ\sm \cU^\pr) \to H^*(\cU, \cU\sm \cU^\pr) \stackrel{f^*}{\to} H^*(U, U\sm U^\pr) \to H^*(Q, Q\sm U^\pr)
\]
sends $\ol{\cD}$ to $\ol{D}$, which proves Lemma \ref{rel}.
\end{proof}
This concludes the proof of Theorems \ref{mainthm} and \ref{surfacethm}.

\section{A generating function}
\label{section:generatingfun}
As was noted in \cite{EGL} and elsewhere, the existence of universal polynomials can be strengthened to a statement about the form of the generating function of Chern integrals.

Throughout this section, we let $X$ be as in the main theorem and let $Q_1, \ldots, Q_k$ be distinct closed, irreducible punctual geometric subsets with $Q_i \subseteq \Hilb_0^{n_i}\left(\CC^d\right)$. 
For any sequence of integers $m_1, \ldots, m_k \ge 0$, let $n_{(m_i)} = \sum m_in_i$.
We let
\[
S_{(m_i)} = \{(j,l)\ |\ 1 \le j \le k, 1 \le k \le m_j\} \subset \ZZ^2.
\]

\begin{ndefn}
We define the geometric set $P((Q_i^{m_i})) \subset X^{[n_{(m_i)}]}$ as the set of $Z = \bigsqcup_{(j,l) \in S_{(m_i)}} Z_{(j,l)}$ such that every $Z_{(j,l)}$ is isomorphic to an element of $Q_j$.
\end{ndefn}
In other words, $P\left((Q_i^{m_i})\right)$ is the set of $Z$ which are the disjoint union of $m_1$ subschemes from $Q_1$, $m_2$ subschemes from $Q_2$, and so on. Specifying appropriate additional data, we can define similar geometric subsets of $X^{\dbl{n_{(m_i)}\dbr}}$.
\begin{ndefn}
\label{def:OrderedGeometricForGenerating}
Let $A = (A_{(j,l)})_{(j,l) \in S_{(m_i)}}$ be a collection of subsets of $\{1, \ldots, n_{(m_i)}\}$, such that $|A_{(j,l)}| = n_j$ and such that the $A_{(j,l)}$ define a partition of $\{1, \ldots, n_{(m_i)}\}$.

We define the geometric set $T((Q_i^{m_i}); A) \subseteq X^{[\![{n_{(m_i)}}]\!]}$ to be the set of all $(Z, (x_i)_{i=1}^{n_{(m_i)}})$ such that $Z = \bigsqcup_{(j,l) \in S} Z_{(j,l)}$, where every $Z_{(j,l)}$ is isomorphic to an element of $Q_j$, and such that $x_i = \Supp Z_{(j,l)}$ if $i \in A_{(j,l)}$.
\end{ndefn}

\begin{nlemma}
\label{thm:DegOfTOverP}
Choosing an $A$ as above, the generic fibre of $T((Q_i^{m_i});A) \to P((Q_i^{m_i}))$ has cardinality $\prod_i m_i!$.
\end{nlemma}
\begin{proof}
We may assume the $Q_i$ are ordered in such a way that if $i < j$, then $Q_i \not\subset Q_j$.
Fix a generic point $Z \in P((Q_i^{m_i}))$.
The number of points in $T((Q_i^{m_i});A)$ lying above $Z$ equals the number of ways of labelling the components of $Z$ by $S_{(m_i)}$ in such a way that $Z_{(j,l)}$ is isomorphic to an element of $S_j$.

Fix one such labelling $Z = \sqcup_{(j,l) \in S_j} Z_{(j,l)}$.
We claim that if $Z = \sqcup_{(j,l) \in S_j} Z^\prime_{(j,l)}$ is a different labelling, then if $Z_{(j,l)} = Z^\prime_{(j^\prime, l^\prime)}$, we must have $j = j^\prime$.

Assume for a contradiction that this is not the case.
Then there must be an equality $Z_{(j,l)} = Z^\prime_{(j^\prime, l^\prime)}$ such that $j < j^\prime$.
We know that $Z_{(j,l)}$ is isomorphic to an element of $Q_j$, and by the ordering of the $Q_j$ that $Q_j \not\subset Q_{j^\prime}$.
Since $Q_j, Q_{j^\prime}$ are closed and irreducible and $Z$ is generic in $P((Q_i^{m_i}))$, it follows that $Z_{j,l}$ is not isomorphic to an element of $Q_{j^\prime}$, which is a contradiction, since $Z_{(j,l)} = Z^\prime_{(j^\prime, l^\prime)} \in Q_{j^{\prime}}$.

Thus the permissible labellings of the components of $Z$ are given by permutations of $S_{(m_i)}$ such that each $(j,l)$ is sent to some $(j,l^\prime)$, of which there are $\prod_i m_i!$.
\end{proof}

Recall that $\CM\left(2,d\right)$ denotes the set of $2$-variable Chern monomials of weight $d$. For a $d$-dimensional $X$ with a bundle $E$ on it, we denote by $N(X,E)$ the Chern number corresponding to an $N \in CM(2,d)$.
\begin{nprop}
\label{genFunProp}
Let $X$ be a smooth, connected projective variety of dimension $d$, let $E$ be an algebraic vector bundle and let $Q_1, \ldots, Q_k$ be distinct closed, irreducible punctual geometric subsets.

We then have
\begin{align*}
\sum_{(m_i) \in \ZZ_{\ge 0}^k} \deg c_\bullet\left(E^{[n_{(m_i)}]}\right) \cap [\overline{P((Q_i^{m_i}))}] x_1^{m_1}\cdots x_k^{m_k} = \prod_{N\in\CM\left(2,d\right)} B_N^{N\left(X,E\right)},
\end{align*}
where the $B_N$ are elements of $\QQ[[x_1, \ldots, x_k]]^\times$ which depend only on the $Q_i$ and the rank of $E$.
\end{nprop}
Note that for degree reasons, the only non-zero term from $c_\bullet\left(E^{[n_{(m_i)}]}\right)$ appearing in the formula is $c_{\dim P((Q_i^{m_i}))}\left(E^{[n_{(m_i)}]}\right)$.

\begin{proof}
Let $F$ be the generating function in the proposition. 
By the main theorem, the coefficients of $F$, and therefore those of $\log F$, are polynomials in the Chern numbers of $(X,E)$.
We will show that the coefficients of $\log F$ are in fact \emph{linear} polynomials in the Chern numbers.
This means that $\log F = \sum N(X,E) b_N$ for $b_N \in \QQ[[y_1, \ldots y_k]]$, and taking $B_N = \exp(b_N)$ will then give the proposition.

Given a sequence $(m_i)_{i = 1}^k$, let $T_{(m_i)} = T((Q_i^{m_i}; A_{(m_i)}) \subset X^{\dbl n_{(m_i)}\dbr}$, where $A_{(m_i)}$ is some appropriate $S_{(m_i)}$-indexed partition of $\{1, \ldots, n_{(m_i)}\}$.
Let
\[
G = \sum_{(m_i) \in \ZZ_{\ge 0}^k} \deg c_\bullet\left(E^{\dbl n_{(m_i)}\dbr}\right) \cap [\overline{T_{(m_i)}}] x_1^{m_1}\cdots x_k^{m_k}.
\]

We use the notation $(m_i)! = \prod_i m_i!$ and denote the coefficient of the $x_{1}^{m_1}\cdots x_{k}^{m_k}$-term of a series by a lower index $(m_i)$.
By Lemma \ref{thm:DegOfTOverP} and the projection formula, we have 
\begin{equation}
\label{eqn:GInTermsOfF}
G_{(m_i)} = (m_i)! F_{(m_i)}.
\end{equation}

Fix now a sequence $(m_i)$, and let $S = S_{(m_i)}$, $n = n_{(m_i)}$ and $A = A_{(m_i)}$.
Let $\mu$ be the partition of $\{1, \ldots, n\}$ induced by $A$, and let $\alpha$ be a partition of $\{1, \ldots, n\}$ such that $\alpha \ge \mu$.
As explained in Section \ref{subsection:approximatingspaces}, for such an $\alpha$ there is an associated $T_\alpha \subset X^{\dbl \alpha \dbr}$ birational to $T$, and we let
\[
G_\alpha = \deg c_\bullet\left(E^{\dbl \alpha\dbr}\right) \cap [\overline{T_\alpha}].
\]

For $1 \le j \le k$, let $S_{j} = \{(j, l)\ |\ 1 \le l \le m_j\} \subset S$.
Giving a partition $\alpha \ge \mu$ is equivalent to giving a partition of $S$, and we will denote this partition of $S$ by $\ol{\alpha}$.
Thus, given a partition $\alpha \ge \mu$, for every $B \in \ol{\alpha}$ we get a sequence $(|B \cap S_j|)_{j = 1}^k$.

As explained in the proof of Lemma \ref{product}, we can decompose $T_\alpha$ and $E^{\dbl \alpha \dbr}$ as products, and in this case the decompositions can be written as $T_\alpha = \prod_{B \in \ol{\alpha}} T_{(|B \cap S_j|)}$ and $E^{\dbl \alpha \dbr} = \bigoplus_{B \in \ol{\alpha}}E^{\left [\!\left[ n_{(|B \cap S_j|)}\right]\!\right]}$.
As in the proof of Lemma \ref{product}, we can now use the Whitney sum formula and the Künneth decomposition to get
\[
\deg c_\bullet\left(E^{\dbl \alpha\dbr}\right) \cap [\ol{T_\alpha}] = \prod_{B \in \ol{\alpha}}\deg c_\bullet\left(E^{\left [\!\left[n_{(|B \cap S_j|)}\right]\!\right]}\right) \cap [T_{(|B \cap S_j|)}]
\]
which implies that 
\begin{equation}
\label{eqn:GAlphaInTermsOfF}
G_\alpha = \prod_{B\in \alpha} G_{(|B \cap S_i|)} = \prod_{B\in \alpha} (|B \cap S_i|)! F_{(|B \cap S_i|)}.
\end{equation}

Let $C_\alpha = c_\bullet(E^{\dbl \alpha \dbr})$, and let $D$ be defined (on the space $Q$ as in the proof of the main theorem) inductively in terms of the $C_\alpha$ as in Definition \ref{def:d}.
By induction one can show that $D = \sum_{\alpha \ge \mu} (-1)^{|\alpha| - 1}(|\alpha|-1)! C_{\alpha}$.
By the argument in Section \ref{sec:relatingClasses}, $\deg D \cap [W]$ is a linear polynomial in Chern numbers of $(X,E)$.
We have
\begin{equation*}
\deg D \cap [W] = \sum_{\alpha \ge \mu} (-1)^{|\alpha| - 1}(|\alpha|-1)! C_\alpha \cap [W] = \sum_{\alpha \ge \mu} (-1)^{|\alpha| - 1}(|\alpha|-1)! G_\alpha.
\end{equation*}
and therefore 
\begin{equation}
\label{eqn:HInTermsOfG}
H_{(m_i)} :=  \sum_{\alpha \ge \mu} (-1)^{|\alpha| - 1}(|\alpha|-1)! G_\alpha
\end{equation}
is a linear combination of Chern numbers.

Using equations \eqref{eqn:GInTermsOfF}, \eqref{eqn:GAlphaInTermsOfF} and \eqref{eqn:HInTermsOfG}, we see that the terms of $H_{(m_i)}$ admits a combinatorial expression in the terms of $F$, and we claim that
\begin{equation}
\label{eqn:HEqualsF}
H_{(m_i)} = (m_i)!(\log F)_{(m_i)}.
\end{equation}
It follows from this that the terms of $\log F$ are linear polynomials in Chern numbers, and so the proposition follows.

We give the proof of \eqref{eqn:HEqualsF} in the case when $k = 1$; the general case can be treated with similar combinatorics.
Fix an $m \ge 0$, and consider $H_m = \sum_{\alpha} (-1)^{|\alpha|-1} (|\alpha| - 1)! \prod_{B \in \alpha}|B|! F_{|B|}$, where the sum is over all partitions $\alpha$ of $\{1, \ldots, m\}$.
Given a partition $\alpha$ of $\{1, \ldots, m\}$, let $\underline{\alpha}$ be the underlying partition of $m$, i.e.\ the sum $\sum_{B \in \alpha} |B| = m$.
For a partition $P = \sum_{i = 1}^l k_i$ of $m$, we use the notation $|P| = l$, $P! = \prod_i k_i!$, $\Aut(P) = \prod_{j \ge 1} |\{i\ |\ k_i = j\}|!$ and $F_P = \prod_i F_{k_i}$. 

Let $P$ be a partition of $m$, and let $H_P = \sum_{\{\alpha | \underline{\alpha} = P\}} (-1)^{|\alpha|-1}(|\alpha| - 1)! \prod_{B \in \alpha} |B|! F_{|B|}$.
All terms in this sum are equal, so we get
\[
H_P = |\{\alpha\ |\ \underline{\alpha} = P\}| \cdot (-1)^{|P|-1}(|P| - 1)! P! F_{P}
\]
The first term of this product equals
\[
\frac{m!}{P! \cdot \Aut(P)}.
\]
We thus get
\[
H_P = \frac{m!}{\Aut(P)}\cdot (-1)^{|P|-1}|P|!/|P| F_{P}.
\]

On the other hand, we can write
\[
\log F = (F-1) - (F-1)^2/2 + \cdots = \sum_{m \ge 0}\sum_{P\vdash m} r_P F_{P} x^{m}
\]
for some integers $r_P$.
All the contributions to $r_P$ come from the term $(-1)^{|P|-1}(F-1)^{|P|}/|P|$.
Expanding out we see that $r_P$ equals $(-1)^{|P|-1}/|P|$ times the number of distinct ways of ordering the terms in $P$, which is $|P|!/\Aut(P)$.

Comparing the terms we see that $H_P = m! r_P F_P$. As $H_m = \sum_{P \vdash m} H_P$ and $(\log F)_m = \sum_{P \vdash m} r_PF_P$, equation \eqref{eqn:HEqualsF} follows.
\end{proof}

\section{Enumerative applications}
\label{applications}
We present some applications of the main theorem to the problem of counting geometric objects with prescribed singularities. We treat three different problems. In Section \ref{section:curvecount} we study curves on a surface having prescribed singularity type, where by singularity type we mean either analytic or topological (equisingular) type. If $L$ is a sufficiently ample line bundle on a surface $S$, we show that the number of such curves in a general linear system $\PP^d \subset |L|$ of appropriate dimension is given by a universal polynomial in the Chern numbers of $(S,L)$. A somewhat more general result to this effect has recently been obtained independently by Li and Tzeng in \cite{LT}.

In Section \ref{divisorcount}, we consider divisors having fixed isolated analytic singularity types on a smooth variety $X$ of arbitrary dimension. We show that the number of such in a general linear system $\PP^d \subset |L|$ is universal. The proofs carry over from the analytic curve singularity case and are omitted.

Finally, in section \ref{BPSsection}, we consider again the case of curves on a surface. We study the locus $|L|_m \subset |L|$ of curves having prescribed ``BPS spectrum'' $m$ and show that if $L$ is sufficiently ample, the Euler characteristic of $|L|_m \cap \PP^k$ is universal.

All of our results use the assumption that $L$ is sufficiently ample. This is required to ensure that the objects we consider occur in the expected codimension in $|L|$, as well as in other places in the argument. A natural way of measuring the ampleness of a line bundle in this setting is $N$-very ampleness, defined as follows.
\begin{ndefn}
 Let $X$ be a nonsingular, projective variety, and let $L$ be a line bundle on $X$. We say that $L$ is $N$-very ample if for every length $(N + 1)$-subscheme $Z \subset X$, the map $H^0(X,L) \to H^0(Z, L|_Z)$ is surjective.
\end{ndefn}	
Equivalently, $L$ is $N$-very ample if the sheaf homomorphism $H^0(X,L) \otimes \cO_{X^{[N+1]}} \to L^{[N+1]}$ is surjective.
Being $0$-very ample is the same as being globally generated, and being $1$-very ample is the same as being very ample.
If $L$ and $M$ are $l$- and $m$-very ample, then $L \otimes M$ is $(l + m)$-very ample \cite{hinohara_tensor_2005}.

\subsection{Curves with specified singularities}
\label{section:curvecount}
We begin by fixing some terms. By a curve singularity we mean a pair $(C,p)$ where $C$ is a reduced, locally planar algebraic curve and $p$ is a singular point of $C$.

Let $\left(C,p\right)$ be a curve singularity. By the analytic type of the singularity $(C,p)$ we mean the isomorphism type of the complete local $\CC$-algebra $\wh{\cO}_{C,p}$. By the topological type (equivalently, equisingularity type) of $\left(C,p\right)$ with $C$ embedded in a smooth surface $S$, we mean the homeomorphism type of the pair $\left(B_\epsilon(p), C \cap B_\epsilon(p)\right)$, where $B_\epsilon(p)$ is a sufficiently small open ball in $S$ centred at $p$.

\begin{nprop}
\label{curveCountProp}
\curveCountProp
\end{nprop}
A similar proposition has recently been obtained by Li and Tzeng in \cite{LT}. They also treat the case where the types are such that some are analytic and some are topological.

\begin{remark}
We will not be concerned with the precise ampleness condition required on $L$ for the universal polynomial to give the correct answer, and refrain from making the $N$ in the statement of the proposition explicit. 
Instead, we will take $N$ large enough whenever the $N$-very ampleness of $L$ is required; it will be clear that $N$ depends only on the types $T_i$. 
A value for $N$ may be recovered from the proof, but already in the case of nodal singularities, the $N$ provided by this method is known to be larger than required by a factor of 5; compare the bounds obtained by our model \cite{gottsche} with those obtained by \cite{KST}.
\end{remark}

The main idea of the proof, taken from \cite{gottsche}, is to set up a correspondence between curves having given singularities and curves containing 0-dimensional subschemes of given isomorphism type. 

Choosing analytic or topological singularity types $T_i$, we show that we may identify curves with such singularities by length $m$ subschemes. We define a geometric set $W = W((T_i)) \subset S^{[m]}$ such that a generic curve containing a $Z \in W$ has the specified singularities. Then, using a proposition from \cite{gottsche} we get that in a general $\PP^d \subseteq |L|$, the number of curves containing a subscheme $Z \in W$ equals $\deg (c_{\dim W}(L^{[m]}) \cap [\ol{W}])$, which is universal by Theorem \ref{mainthm}. We then show that there is a bijection between such pairs $(Z,C)$ and curves in $\PP^d$ with singularities of types $T_i$, completing the proof.

\begin{ncor}
\label{cor:curveGen}
Let $T_1, \ldots, T_k$ be distinct analytic singularity types, and let $n_1, \ldots, n_k$ be non-negative integers. Denote by $G(n_1, \ldots, n_k)$ the universal polynomial computing the number of curves having precisely $n_i$ singularities of type $T_i$ and no other singularities. There are then power series $B_1, B_2, B_3, B_4 \in \QQ[\![x_1, \ldots x_k]\!]$, such that
\[
 \sum G(n_1, \ldots n_k) x_1^{n_1}\cdots x_k^{n_k} = B_1^{c_1^2(L)} B_2^{c_1(L)c_1(S)} B_3^{c_1^2(S)} B_4^{c_2(S)}.
\]
The same statement holds when the $T_i$ are topological types.
\end{ncor}

\begin{proof}
 This follows from Proposition \ref{genFunProp} and the proof of Proposition \ref{curveCountProp}, using the fact that $W$ is irreducible in both the analytic and the topological case.
\end{proof}

\subsubsection{Analytic types}
\label{sec:analyticTypes}
We treat first the case of analytic singularity types. Fix a smooth, projective, connected surface $S$, a line bundle $L$ on $S$, and analytic singularity types $T_1, \ldots, T_k$. We assume that $L$ is $N$-very ample, where $N$ will be taken to be sufficiently large at various points in the proof.

In order to associate a 0-dimensional subscheme to an analytic singularity type, we need the following lemma, which states that a singularity $(C,p)$ is of analytic type $T$ if it looks like a singularity of type $T$ to $M$-th order, where $M$ depends only on $T$.
\begin{nlemma}
Let $(C,p)$ be a curve singularity of analytic type $T$. There is a positive integer $M$, depending only on $T$, such that if $(C^\pr, p^\pr)$ is a curve singularity, the analytic type of $(C^\pr, p^\pr)$ is $T$ if and only if $\cO_{C,p}/\fr{m}^M \cong \cO_{C^\pr, p^\pr}/\fr{m}^M$.
\label{finiteDetermination}
\end{nlemma}
\begin{proof}
This follows from \cite[Cor. 2.24]{GLS} -- in fact we can take $M = \tau +2$, where $\tau$ is the Tjurina number of $T$.
\end{proof}

Given a singularity type $T$ we define the punctual geometric subscheme $W(T) \subset \Hilb^k_0(\CC^2)$ as follows. Let $(C,p)$ be a germ of type $T$, and let $M$ be the integer given by Lemma \ref{finiteDetermination}. Suppose the length of $\cO_{C,p}/\fr{m}_p^M$ is $m$. Let $W(T) \subset \Hilb^m_0(\CC^2)$ be the set of subschemes $Z \in \Hilb_0^m(\CC^2)$ with $Z \cong \Spec \cO_{C,p}/\fr{m}^M$.

Let $m_i$ be the integer such that $W(T_i) \subset \Hilb_0^{m_i}(\CC^2)$, for $i = 1, \ldots, k$. Set $m = \sum m_i$, and define $W \subseteq S^{[m]}$ to be the set of subschemes of the form $Z_1 \sqcup \cdots \sqcup Z_k$, where $Z_i$ is isomorphic to a point in $W(T_i)$ for every $i$.

It is clear that $W$ is a geometric subset, and in the notation of Section \ref{section:geometric} we have $W = P((W(T_i)))$. 
We define the expected codimension of the singularity $T_{i}$ to be $d_{i}= m_{i} - \dim W(T_{i})$.
We let $d = m-\dim W = \sum d_{i}$.

Note that $W(T_i)$ is irreducible and locally closed, as it is the orbit of a given point in $\Hilb^{m_i}_0(\CC^2)$ under the action of the connected algebraic group $\Aut(\cO_{\CC^2,0}/\fr{m}^M_i)$. It follows that $W$ is irreducible and locally closed.

\begin{nlemma}
\label{integralIsIncidenceLemma}
Let $Y \subset S^{[m]}$ be a locally closed subset, and assume $L$ is $(m-1)$-very ample.
\begin{enumerate}
\item [(i)] Let $\cZ \subset S^{[m]} \times |L|$ denote the incidence locus of pairs $(Z,C)$ with $Z \in Y$ and $Z \subset C$. We have $\dim \cZ = \dim |L| + \dim Y - m$.
\item[(ii)] Let $e= m - \dim Y$, and let $\PP^e \subset |L|$ be a general linear subspace. The number of pairs $(Z,C)$ such that $Z \in Y$, $C \in \PP^e$ and $Z \subset C$ is equal to
\[
 \deg c_{\dim Y}(L^{[m]}) \cap [\ol{Y}]
\]
\end{enumerate}
\end{nlemma}

\begin{proof}
(i) For any $Z \in Y$, the fibre of $\cZ \to Y$ over $Z$ is the projectivisation of the kernel of $H^0(S,L) \to H^0(Z,L|_Z)$. By the $(m-1)$-very ampleness of $L$, this homomorphism is surjective, so $\cZ \to Y$ is a projective space bundle with fibres of dimension $|L|-m$. The claim follows.

(ii) See the proof of \cite[Prop 5.2]{gottsche}.
\end{proof}

Applying Lemma \ref{integralIsIncidenceLemma} (ii) with $W = Y$, the following lemma now concludes the proof of Proposition \ref{curveCountProp} in the analytic case.
\begin{nlemma}
\label{lemma:analyticbijection}
Let $\PP^d \subset |L|$ be a general subsystem, and assume $L$ is $N$-very ample. Suppose $(Z,C)$ is a pair such that $Z \in W$, $C \in \PP^d$ and $Z \subset C$. Then $C$ has $k$ singularities of analytic types $T_i$, and $C$ contains no other point of $W$.
\end{nlemma}
\begin{proof}
\newcommand{\ord}{\text{ord}}
Let $\PP^d$, $C$ and $Z$ be as in the statement of the lemma, and suppose $Z = \sqcup Z_i$, where $Z_i$ is supported at $x_i \in C$ and where $Z_i$ is isomorphic to a point in $W(T_i)$. We show the following claims: (1) That $C$ has precisely $k$ singularities, (2) that $C$ has a singularity of type $T_i$ at $x_i$, and (3) that $C$ contains precisely one $Z \in W$.

(1): Clearly, $C$ has at least $k$ singularities. Assume for a contradiction that $C$ has more than $k$ singularities. It must then contain a subscheme of the form $Z \sqcup Z^\pr$, where $Z^\pr$ is defined by an ideal $\fr{m}_x^2$ for some $x \in S$ where $C$ is singular. The geometric set
\[
 W^\pr := \{Z \sqcup Z^\pr \mid Z \in W \text{ and } Z^\pr = \Spec\cO_{S,x}/\fr{m}_x^2\} \subset S^{[m+3]}
\]
has dimension 2 greater than $W$. By Lemma \ref{integralIsIncidenceLemma} (i), we see that the set of $C \in |L|$ containing an element of $W^\pr$ has codimension $>d+1$ in $|L|$ if $L$ is $(m+2)$-very ample. The intersection of this set with a general $\PP^d \subset |L|$ is empty, contradicting the original assumption.

(2): Suppose for a contradiction that the singularity type of $C$ at $x_1$ is $T_1^\pr \not= T_1$. Let $M$ be the integer associated to $T_1$ as in Lemma \ref{finiteDetermination}, and let $R =\cO_{S,x_1}/\fr{m}^M$. As $T_1^\pr \not= T_1$, we have $Z_1 \subsetneq C\cap \Spec R$. Let $f,g \in R$ be defining equations of $Z_1$ and $C \cap \Spec R$ in $R$, we then have $(g) \subsetneq (f)$. This implies that $(g) \subseteq \fr{m}\cdot(f) \subseteq \fr{m}^{\ord(f)+1} \cap (f)$, where $\text{ord}(f)$ is the maximal integer such that $f \in \fr{m}^{\text{ord}(f)}$. 

Hence $C$ contains a subscheme of the form
\[
 Z_1^\pr \sqcup Z_2 \sqcup \ldots \sqcup Z_k,
\]
where $Z_1^\pr = \Spec R/\fr{m}^{\ord(f)+1} \cup Z_1$. Let $W^\pr$ be the set of subschemes which can be written in this way. Then $W^\pr$ is geometric and has dimension $\le \dim W$.

Let $m^\pr$ be the length of the points of $W^\pr$. Clearly, we have $m^\pr > m$, so we have
\[
m^\pr - \dim W^\pr > m -\dim W = d.
\]
As $L$ is $N$-very ample, applying Lemma \ref{integralIsIncidenceLemma} (i) shows that the codimension of the locus of points containing a point from $W^\pr$ is $> d$, if $N \ge m + 1$. As $C \in \PP^d$ for a general $\PP^d$, it is not contained in this locus.

(3): By (1) and (2), we know that $C$ has a singularity of type $T_i$ at $x_i$, and suppose for a contradiction that there is a $Z^\pr \in W$ with $Z^\pr \subset C$ such that $Z^\pr \not= Z$. Let $Z^\pr = \sqcup Z^\pr_i$ with $Z^\pr_i$ supported at $x_i$. Assume that $Z_i \not = Z_i^\pr$ as subschemes. But by part (2), the singularity type associated to $Z_i^\pr$ must be $T_i$, and hence we have
\[
 Z_i = \Spec \cO_{C,x_i}/\fr{m}^M = Z^\pr_i,
\]
where $M$ is as in Lemma \ref{finiteDetermination} for type $T_i$.
\end{proof}

\subsubsection{Topological singularities}
\label{subsubsection:topSing}
We now turn to the case of topological singularities. 
Let $S$ and $L$ be as before, and fix topological singularity types $T_1, \ldots, T_k$.
For any planar curve singularity $(C,p)$ the infinitely near points in $C$ over $p$ define a combinatorial structure called the Enriques diagram, which determines the equisingularity type of $(C,p)$ \cite{KP}.
Let $D$ be the Enriques diagram of the $T_{i}$, that is the union of the Enriques diagrams for each $T_{i}$.

The degree of an Enriques diagram is defined in \cite{KP}, and we let $m = \deg(D)$.
Kleiman and Piene construct a subscheme $H(D) \subset S^{[m]}$.
It has the property that if $Z \in H(D)$ and $C$ is a generic curve containing $Z$, then the singularities of $C$ correspond to the Enriques diagram $D$.
The subset $H(D)$ is geometric and irreducible \cite[5.8]{KPT}.

Let $d = m - \dim H(D)$.
The following lemma is a reformulation of \cite[3.7]{KP}.
\begin{nlemma}
Let $\PP^d \subset |L|$ be a general subsystem, and assume $L$ is $(m-1)$-very ample. Suppose $(Z,C)$ is a pair such that $Z \in H(D)$, $C \in \PP^d$ and $Z \subset C$. Then $C$ has $k$ singularities of topological types $T_i$, and $C$ contains no other point of $H(D)$.
\end{nlemma}
\begin{proof}
Let $Y$ be the incidence locus in $|L| \times W$, let $|L|_{T} \subset |L|$ be the set of curves having prescribed singularity types, let $\pi : |L| \times W \to |L|$ be the projection and let $Y_{T} = \pi^{-1}(|L|_{T}) \subset Y$.
By \cite[3.7]{KP}, $Y_{T}$ is dense in $Y$ when $L$ is $(m-1)$-very ample.\footnote{The reference assumes $L = M \otimes N^{\otimes m}$ for $M$ globally generated and $N$ very ample. However, the proof given shows that $L$ is then $(m-1)$-very ample, and the stronger assumption on $L$ is not needed.}
We therefore have $\dim Y \backslash Y_{T} < \dim Y$, and applying Lemma \ref{integralIsIncidenceLemma} we find $\dim Y = \dim |L| - d$.
So if $\PP^{d} \subset |L|$ is general, then $\pi^{-1}(\PP^{d}) \cap Y \subset Y_{T}$, proving the first claim of the lemma.
By \cite[3.7]{KP}, the map $Y_{T} \stackrel{\pi}{\to} |L|_{T}$ is bijective, which proves the second claim.
\end{proof}

Applying Lemma \ref{integralIsIncidenceLemma} (ii) with $Y = H(D)$, the above lemma concludes the proof of Proposition \ref{curveCountProp} in the topological case.

If $D$ is a connected Enriques diagram, then every element $Z \in H(D)$ will have support in one point.
Since $H(D)$ is geometric, it is defined by a punctual geometric subset $H(D)_{0} \subset \Hilb^{\deg D}_{0}(\CC^{2})$.

For topological singularites $T_{1}, \ldots, T_{k}$, the associated Enriques diagram is $D = D_{1} \sqcup \ldots \sqcup D_{k}$, where $D_{i}$ is the connected Enriques diagram associated with $T_{i}$.
In the notation of Section \ref{section:geometric}, we have $H(D) = P((H(D_{i})_{0}))$, so we can apply Proposition \ref{genFunProp} to prove Corollary \ref{cor:curveGen} for the case of topological singularity types.

\subsection{General hypersurface singularities}
\label{divisorcount}
Without any extra work, the above extends to counts of analytic types of isolated singularities of hypersurfaces. Let $(D,x)$ be the pair of a divisor $D$ in a nonsingular variety $X$ and an isolated singular point of $D$. 
The analytic type of the singularity $(D,p)$ is the isomorphism type of the complete local $\CC$-algebra $\wh{\cO}_{D,x}$.

Lemma \ref{finiteDetermination} is valid for hypersurface singularities of all dimensions, and the following proposition can be shown by the proof given in Section \ref{sec:analyticTypes}, mutatis mutandis.
(Note in particular that we do not use the nonsingularity of $X^{[n]}$ anywhere in the argument.)
\begin{nprop}
\label{divCountProp}
\divCountProp
\end{nprop}

There is a similar corollary for the generating function of these universal polynomials as in the curve case; we leave the statement of this to the reader.

\subsection{BPS spectrum loci}
\label{BPSsection}
Let $C$ be a reduced, complete, locally planar algebraic curve, and consider the generating function
\[
 H_C(q) = \sum_{k=0}^\infty \chi\left(C^{[k]}\right) q^k.
\]
Let the arithmetic and geometric genus of $C$ be $g(C)$ and $\ol{g}(C)$, respectively. In \cite{PT} it is shown that there are integers $n_{i,C}$, with $n_{i,C} = 0$ unless $\ol{g} \le i \le g$, such that
\begin{equation}
\label{equation:BPS}
 H_C(q) = \sum_{i=\ol{g}(C)}^{g(C)} n_{i,C} q^{g-i}(1-q)^{2i-2}.
\end{equation}

For our purposes, it will be convenient to work with the index-shifted integers $m_{i,C} := n_{g-i,C}$. We define the BPS spectrum of $C$ to be the sequence of integers $(m_{i,C})_{i=0}^\infty$. 
By the above, we have $m_{i,C} = 0$ if $i \ge g -\ol{g}$. 
If $C$ has $k$ singularities of analytic types $T_1, \ldots, T_k$, then by stratifying $C^{[k]}$ one can see that the BPS spectrum of $C$ depends only on the $T_i$.

By this observation, one may define the BPS spectrum of an analytic singularity type $T$ as the BPS spectrum of a complete, reduced curve having one singularity of type $T$. The BPS spectrum of a singularity $T$ is shown by Maulik \cite{maulik_stable_2012} to be determined explicitly by the Milnor number and the HOMFLY polynomial of the link of $T$. In particular, the BPS spectrum of a curve depends only on the topological types of the singularities of the curve.

\begin{nprop}
 \BPSProp
\label{prop:BPS}
\end{nprop}
If in addition it is known that $\PP^k \cap |L|_m$ is 0-dimensional, this implies an enumerative result of the kind found in the previous subsection. This is essentially the argument used in \cite{KST} to compute the number of $\delta$-nodal curves and prove the Göttsche Conjecture.

The remainder of this section contains the proof of Proposition \ref{prop:BPS}.
\begin{nlemma}
\label{thm:FinitelyManyValuesOfChi}
Let $j \ge 0$ be an integer. Then $\chi(\Hilb^j_p(C))$ can take only finitely many values for $(C,p)$ a locally planar curve singularity.
\end{nlemma}
\begin{proof}
Let $R = \CC[x,y]$ and let $R_{<j}$ denote the vector space of polynomials of degree $< j$.
Let $B = \Spec\CC[x_m]_{m \in R_{<j}}$, and let $\cC \subset \AA^2 \times B$ be the divisor defined by $\cI = \sum_{m \in R_{<j}} m x_m$.
Let $\Hilb^j_0(\cC/B)$ be the relative punctual Hilbert scheme parametrising length $j$ subschemes which are contained in some fibre $\cC_b$ of $\cC \to B$ and are then supported at the point $(0,b) \in \cC_b$.

There is a map $f:\Hilb^j_0(\cC/B) \to B$.
For a given planar curve singularity $(C,p)$, the scheme $\Hilb^j_p(C)$ depends only on the $(j-1)$-st order neigbourhood $\Spec(\cO_{C,p}/\mathfrak{m}_p^j)$.
It follows that $\Hilb^j_p(C)$ is isomorphic to some fibre of $f$.
The Euler characteristic push-forward function $f_*(1)$ is a constructible function on $B$, hence takes only finitely many values.
The claim follows. 
\end{proof}

\begin{nlemma}
\label{thm:genusCurve}
If $L$ is $k$-very ample, then for a general $\PP^k \in |L|$ every curve $C \in \PP^k$ is reduced and satisfies $\ol{g}(C) \ge g(C) - k$.
\end{nlemma}
\begin{proof}
See \cite[Prop.\ 2.1]{KST}.
\end{proof}

\begin{nlemma}
Let $k$ be an integer. There is a finite set $T$ of BPS spectra such that if $L$ is $k$-very ample and $\PP^k \subset |L|$ is a general linear subsystem, then for every curve $C \in\PP^k$ the BPS spectrum of $C$ is in $T$.
\end{nlemma}
\begin{proof}
By Lemma \ref{thm:genusCurve}, we have $\ol{g}(C)\ge g(C) - k$. By \eqref{equation:BPS}, the BPS spectrum of $C$ is then determined by $\chi(C^{[i]})$ for $1 \le i \le k$. 

Denote by $\Hilb^j_p(C) \subset C^{[j]}$ the set of subschemes supported at $p \in C$. Stratifying $C^{[i]}$, we see that $\chi(C^{[i]})$ is determined by $\chi(C)$ and the integers $\chi(\Hilb^j_p(C))$, where $p$ ranges over the singular points of $C$ and $j \le i$. Applying Lemma \ref{thm:FinitelyManyValuesOfChi}, the claim follows.
\end{proof}

\begin{nlemma}
Let $k, L$ and $T$ be as in the previous lemma, and let $m \in T$ be a BPS spectrum.
There is an $F_m \in \QQ(g)[x_1, \ldots, x_k]$ such that
\[
 F_m\left(g(C), \chi\left(C\right), \ldots, \chi\left(C^{[k]}\right)\right)
\]
equals 1 if $C$ has BPS spectrum $m$ and equals $0$ if $C$ has BPS spectrum in $T \setminus \{m\}$.
\end{nlemma}
\begin{proof}
Let $C$ be a curve with BPS spectrum $m$.
Using \eqref{equation:BPS}, we find that $\chi(C^{[i]})$ is a polynomial in $g(C)$:
\[
\chi(C^{[i]}) = \sum_{j=0}^i(-1)^{i-j}m_{j} \binom{2(g(C)-j)-2}{i-j}.
\]

By Lemma \ref{thm:genusCurve}, we have $\ol{g}(C)\ge g(C) - k$. 
By \eqref{equation:BPS}, the BPS spectrum of $C$ is then determined by $\chi(C^{[i]})$ for $1 \le i \le k$.

For all $m \in T$, let $C_m$ be a curve of BPS spectrum $m$.
Let $P_m = (\chi(C_m^{[1]}), \ldots, \chi(C_m^{[k]})) \in \QQ(g)^k$.
The $P_m$ are all distinct, and for each $m$ we can find an element $G_m \in \QQ(g)[x_1, \ldots, x_k]$ such that $G_m(P_{m^{\prime}}) = 0$ if and only if $i = j$.
Putting $F_i = \prod_{j \not= i} \frac{G_j}{G_j(P_i)}$ gives the result.
\end{proof}

Let now $\PP^k \subset |L|$ be general, let $\cC \to \PP^k$ be the family of curves, and let $\cC^{[i]}/\PP^k$ denote the relative Hilbert scheme. Every monomial $r$ in the variables $x_i$ determines a scheme $\cC\left(r\right)$ by taking
\[
 \cC\left(x_i\right) = \cC^{[i]}/\PP^k
\]
and extending this by the rule
\[
 \cC\left(r_1 \cdot r_2\right) = \cC\left(r_1\right)\times_{\PP^k}\cC\left(r_2\right).
\]
It is clear that
\begin{align*}
 \chi\left(\cC\left(r\right)\right) &= r(\chi(\cC^{[1]}/\PP^k), \ldots, \chi(\cC^{[k]}/\PP^k)) \\
&=\sum_{m\in T} \chi\left(|L|_{m} \cap \PP^k \right) r\left(\chi\left(C_m^{[1]}\right), \ldots, \chi\left(C_m^{[k]}\right)\right),
\end{align*}
where $C_{m}$ denotes a curve with BPS spectrum $m$.

We may write $F_i$ in the form
\[
F_i = \sum_{r} f_r(g) r(x_1,\ldots, x_k)
\]
where the sum is over monomials $r$ and where $f_r \in \QQ(g)$. Since $g(C) = (c_1(L)^2 - c_1(L)c_1(S))/2 + 1$, we get
\[
\chi(|L|_{m} \cap \PP^k) = \sum_r f_r((c_1(L)^2 - c_1(L)c_1(S))/2 + 1) \chi(\cC(r)).
\]
Lemma \ref{thm:BPSIntegralsAreUniversal} below shows that $\chi(\cC(r))$ is universal.
This implies that $\chi(|L|_{m} \cap \PP^k)$ admits a universal expression as $G/H$, where $G$ is a polynomial in the Chern numbers of $(S,L)$ and $H$ is a polynomial in $g = (c_1(L)^2 - c_1(L)c_1(S))/2 + 1$.
In Section \ref{sec:bootstrap} we amplify this and show that in fact $G/H$ is a polynomial, concluding the proof of Proposition \ref{prop:BPS}.

\begin{nlemma}
\label{thm:BPSIntegralsAreUniversal}
Let $L$ be a line bundle, and let $\PP^k \subseteq |L|$ be a general linear subsystem. Let $\cC \to \PP^k$ be the universal family of curves, and denote by $\cC^{[i]} \to \PP^k$ the relative Hilbert scheme of $i$ points for this morphism. Then the Euler characteristic
\[
 \chi\left(\cC^{[i_1]} \times_{\PP^k} \cdots \times_{\PP^k} \cC^{[i_l]}\right)
\]
is computed by a universal polynomial, provided that $L$ is $((\sum_j i_j) - 1)$-very ample.
\end{nlemma}
\begin{proof}
For notational simplicity, we treat the case where $l = 2$, the general case is essentially the same. The case $l = 1$ is simpler, see \cite{KST}. 

Let $f_n : \cC^{[n]} \to S^{[n]}$ be the natural morphism.
We claim that there exists a finite stratification of $S^{[n]}$ by geometric sets $W_{n,i}$ of universal type, such that
\begin{equation}
\label{wRel}
\chi(\cC^{[k_1]}\times_{\PP^k}\cC^{[k_2]}) = \sum_{n=1}^{k_1+k_2}\sum_{i} i\cdot \chi(\cC^{[n]} \cap f_n^{-1}(W_{n,i})).
\end{equation}
Consider the function $g:S^{[k_1]} \times S^{[k_2]} \to \sqcup_{n=1}^{k_1+k_2} S^{[n]}$ defined pointwise by
\[
 g(Z_1,Z_2) = Z_1 \cup Z_2,
\]
where the union is in the scheme-theoretic sense. 

Define
\[
 W_{n,i} = \{Z \in S^{[n]} \mid \chi(g^{-1}(Z)) = i)\}.
\]
One can check that $W_{n,i}$ is geometric. Using the fact that $Z_1, Z_2 \subset C \Leftrightarrow Z_1 \cup Z_2 \subset C$, we also see that $W_{n,i}$ satisfies \eqref{wRel}. Lemma \ref{lemma:universalChi} now completes the proof.
\end{proof}

\begin{nlemma}
\label{lemma:universalChi}
Let $W$ be a geometric subset of $S^{[n]}$, and let $\PP^k \subset |L|$ be a general linear subsystem, with $L$ an $(n-1)$-very ample line bundle. Let $\cC \to \PP^k$ be the family of curves, let $\cC^{[n]}$ be the relative Hilbert scheme of the family, and let $f:\cC^{[n]} \to S^{[n]}$ be the natural morphism. 

Then there exists a universal polynomial in the Chern numbers of $(S,L)$ which computes $\chi\left(f^{-1}(W) \cap \cC^{[n]}\right)$.
\end{nlemma}

\begin{proof}
The inclusion-exclusion principle for $\chi$ lets us reduce to the case where $W$ is closed and irreducible. Consider the diagram
\[
 \begin{CD}
  f^{-1}(W) \cap \cC^{[n]} @>f>> W \\
  @VVV \\
  \PP^k.
 \end{CD}
\]
The fibres of $f$ are all projective spaces, since for a point $Z \in W$, the fibre over $Z$ is the linear system of curves containing $Z$. Hence we have $\chi\left(f^{-1}\left(Z\right)\right) = \dim f^{-1}\left(Z\right) + 1$. Let $W_m = \{Z \mid \chi(f^{-1}(Z)) = m\}$. On $W$, consider the surjective homomorphism
\[
 H^0\left(S,L\right)\otimes \cO_{S^{[n]}} \to L^{[n]},
\]
let $E \subseteq H^0\left(S, L\right)$ be the $\left(k + 1\right)$-dimensional subspace defining $\PP^k$, and let $\phi : E\otimes \cO_W \to L^{[n]}$ be the induced homomorphism. Then 
\[
W_m = \{Z \in W \mid \dim \ker \phi = m\}.
\]
Letting $D_r(\phi)$ denote the locus over which $\phi$ has rank $\le r$, we have $W_m = D_{k+1-m}(\phi)\sm D_{k+1-m-1}(\phi)$. It thus suffices to compute $\chi(D_r(\phi))$ for all $r$.

By \cite[Thm 2.10]{PP}, there exists a formula for the Euler characteristic of $D_r(\phi)$ as a polynomial in the Chern classes of $L^{[n]}$ capped with $\csm(W)$, assuming the homomorphism $E \to L^{[n]}$ is $r$-general in the sense of \cite{PP}. Choosing a Whitney stratification of $W$, $r$-generality amounts to saying that over each stratum, the section of $\sHom\left(E, L^{[n]}\right)$ intersects the tautological degeneracy locus $D_r \setminus D_{r-1} \subseteq \sHom\left(E, L^{[n]}\right)$ transversely.

The $(n-1)$-very ampleness of $L$ implies there is a surjection $H^0\left(S,L\right) \otimes \cO_{W} \to L^{[n]}$, inducing a morphism
\[
 W \to \Gr\left(H^0\left(S,L\right), n\right).
\]
Choosing a subspace $E \subseteq H^0\left(S,L\right)$, the intersection of $\phi$ with $D_r\sm D_{r-1}$ corresponds to the intersection of $W$ with a certain smooth subset of $\Gr\left(H^0\left(S,L\right),n\right)$. By the Kleiman-Bertini transversality theorem \cite{kleiman}, for a general $E \subseteq H^0\left(S,L\right)$, the intersection of each Whitney stratum of $W$ with this set will be smooth of the expected dimension. This shows that $E\otimes\cO_{W} \to L^{[n]}$ is $r$-general in the sense of \cite{PP}.

Hence the formula of \cite[Thm 2.10]{PP} applies, and by Theorem \ref{mainthm} (ii), the statement of the lemma follows.
\end{proof}

\subsection{A bootstrap}
\label{sec:bootstrap}
Let $A = \QQ[x_{c_{1}(S)^{2}}, x_{c_{1}(S)c_{1}(L)}, x_{c_{1}(L)^{2}}, x_{c_2(S)}]$.
For any $F \in A$, write $F(S,L)$ for the value obtained by assigning the Chern numbers of $(S,L)$ to the $x_i$.
We have shown above that there are $G, H \in A$ such that $G(S,L)/H(S,L)$ computes the Euler characteristic of $\chi\left(\PP^k \cap |L|_{m}  \right)$ when $L$ is $N$-very ample.
Furthermore, $H$ is contained in the subring $\QQ[x_{g}]$, where $x_{g} = (x_{c_{1}(L)^{2}}-x_{c_{1}(L)c_{1}(S)})/2 - 1$.

The claim of Proposition \ref{prop:BPS} is that we may take $H=1$ here, or equivalently that $H$ divides $G$.
Lemma \ref{thm:HDividesG} below shows that this is indeed the case.

\begin{nlemma}
\label{thm:RationalFunctionIsPolynomial}
Let $F \in \QQ(x)$ be such that $F(n)$ is an integer for all $n \gg 0$. Then $F \in \QQ[x]$.
\end{nlemma}
\begin{proof}
We can write $F = Q + R/H$ with $H, Q, R \in \QQ[x]$ and with $\deg R < \deg H$.
Let $N \in \ZZ$ be such that $NQ \in \ZZ[x]$; then $NF(n) - NQ(n)$ is integral for $n \gg 0$.
This equals $NR(n)/H(n)$, which tends to 0 as $n \to \infty$.
Hence $R(n) = 0$ for $n \gg 0$, and so $R = 0$.
\end{proof}

\begin{nlemma}
\label{thm:surfacesZariskiDense}
The set of quadruples
\[
(c_{1}(L)^{2}, c_{1}(L)c_{1}(S), c_{1}(S)^{2}, c_{2}(S)) \in \ZZ^{4}
\]
where $S$ is a smooth, connected, projective surface and $L$ is an ample line bundle forms a Zariski dense subset of $\CC^{4}$.
\end{nlemma}

\begin{proof}
Let $\cS$ be the set of surfaces such that the Picard rank is $\ge 2$ and $c_{1}(S)$ is numerically non-trivial.
Let $S \in \cS$ and let $NS(S)$ be the Neron-Severi group of $S$.
Since the ample classes in $NS(S)\otimes \RR$ form an open cone, we see that $\{[L]\ |\ L \text{ ample}\}$ is a Zariski dense subset of $NS(S) \otimes \CC$.
One checks that the set $\{(\alpha^{2}, \alpha c_{1}(S))\ |\ \alpha \in NS(S) \otimes \CC\}$ is Zariski dense in $\CC^{2}$, and hence $\{(c_{1}(L)^{2}, c_{1}(L)c_{1}(S))\ |\ L \text{ ample}\}$ is Zariski dense in $\CC^{2}$.

Therefore the Zariski closure of the set in the statement of the lemma contains 
\[
\{(a,b,c_{1}(S)^{2},c_{2}(S)\ |\ a,b \in \CC, S \in \cS\}.
\]
Within $\cS$ is the set $\cS^{\pr}$ of surfaces birational to the product of two curves of genera $\ge 2$.
One checks that $\{(c_{1}^{2}(S), c_{2}(S))\ |\ S \in \cS^{\pr}\}$ is Zariski dense in $\CC^{2}$.
The claim follows.
\end{proof}

\begin{nlemma}
\label{thm:HDividesG}
Assume $G \in A$ and $H \in \QQ[x_{g}] \subset A$ are such that if $L$ is $N$-very ample, then $G(S,L)/H(S,L)$ is an integer. Then $H$ divides $G$.
\end{nlemma}
\begin{proof}
Consider the map $A \to A[t]$ which sends $x_{c_{1}(L)^{2}} \mapsto t^{2} x_{c_{1}(L)^{2}}$, $x_{c_{1}(L)c_{1}(S)} \mapsto tx_{c_{1}(L)c_{1}(S)}$ and leaves the other generators fixed.
Let $G^{\prime}, H^{\prime} \in A[t]$ be the images of $G,H$ under this map.
Then if $t \in \ZZ$, we have $G^{\prime}(S,L)(t) = G(S,tL)$ and likewise for $H$.

Ordering by $t$-degree, the leading term of $H^{\prime}$ is proportional to $t^{2k}x^{k}_{c_{1}(L)^{2}}$ for some $k$.
We may therefore write $G^{\prime} = QH^{\prime} + R$, where $Q, R \in A[x^{-1}_{c_{1}(L)^{2}}, t]$ and the $t$-degree of $R$ is less than that of $H$.

Let $(S,L)$ be such that $L$ is ample.
For $t \in \ZZ$ we have $G^{\prime}(S,L)(t)/H^{\prime}(S,L)(t) = G(S,tL)/H(S,tL)$.
If $t \gg 0$, then $tL$ is $N$-very ample, and so $G(S,tL)/H(S,tL)$ is an integer.
By Lemma \ref{thm:RationalFunctionIsPolynomial}, $H^{\prime}(S,L)$ must then divide $G^{\prime}(S,L)$.
It follows that $R(S,L) = 0$.
Since this is true for all pairs $(S,L)$ with $L$ ample, \ref{thm:surfacesZariskiDense} implies that $R = 0$.

We thus have $G^{\prime} = QH^{\prime}$.
Setting $t = 1$ gives $G = PH$, where $P = Q|_{t=1} \in A[x^{-1}_{c_{1}(L)^{2}}]$.
Since $H \in \QQ[x_{g}]$, it is not divisible by $x_{c_{1}(L)^{2}}$, and it follows that $P \in A$.
This proves the claim.
\end{proof}

\footnotesize{
\bibliographystyle{alpha-abbrv}
\bibliography{UniversalThm}
}
\end{document}